\documentclass[preprint,review, 12pt]{elsarticle}
\usepackage{graphicx} 
\usepackage{amsmath,amssymb,mathtools}

\usepackage{amssymb}
\usepackage{epsfig}
\usepackage{comment}
\usepackage{amsmath}
\usepackage{subcaption}
\usepackage{blindtext}
\usepackage{xcolor}
\usepackage{hyperref}
\usepackage[section]{placeins}
\usepackage{mathrsfs}
\usepackage{graphicx}
\usepackage[notrig]{physics}
\usepackage{units}
\usepackage{stmaryrd}
\usepackage{algorithm,algorithmic}

\newtheorem{thm}{Theorem}[section]

\newtheorem{example}[thm]{Example}

\DeclareMathOperator{\dist}{dist}

\newcommand{\D}{\mathcal{D}}

\newcommand{\PM}{\mathcal{P}}
\newcommand{\hF}{\hat{F}}

\newtheorem{theorem}{Theorem}

\newcommand{\R}{\mathbb{R}}

\usepackage{todonotes}

\begin{document}

\begin{frontmatter}

\title{Approximating Young Measures With Deep Neural Networks}
\author[inst1,inst2]{Rayehe Karimi Mahabadi}
\affiliation[inst1]{organization={Department of Civil and Environmental Engineering, Duke University},
            city={Durham}, 
            state={NC},
            country={USA}}
            
\affiliation[inst2]{organization={Department of Mechanical Engineering and Materials Science, Duke University},
            city={Durham}, 
            state={NC},
            country={USA}}

\author[inst3,inst4,inst5]{Jianfeng Lu\corref{cor1}}
\ead{jianfeng@math.duke.edu}
\affiliation[inst3]{organization={Department of Mathematics, Duke University},
            city={Durham}, 
            state={NC},
            country={USA}}

\affiliation[inst4]{organization={Department of Physics, Duke University},
            city={Durham}, 
            state={NC},
            country={USA}}

\affiliation[inst5]{organization={Department of Chemistry, Duke University},
            city={Durham}, 
            state={NC},
            country={USA}}

\author[inst1,inst2]{Hossein Salahshoor\corref{cor1}}
\ead{hossein.salahshoor@duke.edu}

\cortext[cor1]{Corresponding authors}

\begin{abstract}
Parametrized measures (or Young measures) enable to reformulate non-convex variational problems  as convex problems at the cost of enlarging the search space from space of functions to space of measures. To benefit from such machinery, we need powerful tools for approximating measures. We develop a deep neural network approximation of Young measures in this paper. The key idea is to write the Young measure as push-forward of Gaussian measures, and reformulate the problem of finding Young measures to finding the corresponding push-forward. We approximate the push-forward map using deep neural networks by encoding the reformulated variational problem in the loss function. After developing the framework, we demonstrate the approach in several numerical examples. We hope this framework and our illustrative computational experiments provide a pathway for approximating Young measures in their wide range of applications from modeling complex microstructure in materials to non-cooperative games. 
\end{abstract}

\end{frontmatter}
\section{Introduction}

From an energy minimization perspective, physical systems often represent several equilibrium states which confers a multi-well structure in their energy landscape. This ubiquitously  emerges in mechanics  and material science and, in turn, renders the resultant variational problems as non-convex. There are many examples of non-convex variational problem in continuum mechanics, from optimal design problems \cite{goodman1986numerical}, to micromagnetics \cite{de1993energy}, and to crystalline materials \cite{luskin1996computation,ball1987fine,chipot1988equilibrium}. A prototypical example is the theory of martensitic microstructure where many equilibrium phases co-exist and crystal structure can alternate between these multiple phases --- the interested reader is referred to \cite{bhattacharya2003microstructure} and references therein. 

From a calculus of variations point of view, direct methods are not applicable to such non-convex problems, as the minimizing sequences may not converge in the strong topology. In practice, the corresponding minimizing sequences often develop ever increasing oscillations between the multiple equilibria. Solutions to such problems often fail to converge with mesh refinement, such that a solution obtained at one resolution may completely change when computed on a finer mesh. Numerical schemes therefore yield coarse, mesh-biased approximations that suppress many competing microstructures, frequently trapping the solution in metastable states that misrepresent materials' true energetics and deformation.

To tackle these non-convex problems, one often derives a modified variational problem through a relaxation procedure \cite{chipot1988equilibrium, kohn1991relaxation, govindjee2001multi, luskin1996computation}. In such approaches, the original problem is replaced by a \textit{(quasi)convex envelop}, that is the lowest energy possible through all possible \textit{microstructure}. Obtaining explicit relations in this quasiconvexification step is only achievable, if at all possible, for a very tiny class of microstructures such as laminates \cite{kohn1991relaxation,de1993energy,luskin1996computation,ortiz1999nonconvex,aubry2003constrained}. Alternatively, one can use Young measures. These parameterized measures are objects that describe the limit of minimizing sequences and allow for directly finding the effective energy without going through the relaxation procedure. 

Young measures are inherently high-dimensional objects and difficult to approximate. Thus, despite their nice analytical properties, numerical approaches based on Young measures are less explored. Initial efforts have been made to directly approximate Young gradient measures as a convex combination of Dirac masses \cite{nicolaides1993computation,carstensen2000numerical}, with the locations of the masses fixed on a uniform grid in phase space and the admissible class represented by finite element spaces. On the other hand, in recent years, deep neural networks have shown strong promise in treating very high-dimensional problems  \cite{raissi2019physics,karniadakis2021physics, han2018solving, beck2019machine,  li2020fourier, kovachki2023neural, khoo2021solving,lu2020universal}. In this paper, we aim to leverage the power of neural networks and develop a deep learning based framework for approximating Young measures.

Here is the road-map for the rest of the paper: in section 2, we re-visit the fundamentals of Young measure theory and how they naturally appear in non-convex variational problems. Section 3 focuses on developing a neural network representation of Young measures. We then apply this framework for solving a number of numerical examples. We conclude by some remarks on potential applications of our framework in mechanics and material science. 

\section{Young measures and non-convex optimization}

Young measures are maps from a domain to space of probability measures. They were introduced in the pioneering works of L. C. Young \cite{young1937generalized,young2024lectures}, and since then have found many applications in 
homogenization theory \cite{ tartar1983compensated,tartar1990h,allaire1994optimal, ball2005version, francfort1994sets, milton1990characterizing, milton1994link, nesi1991polycrystalline, alberti2001new}, optimal control \cite{balder1984general, roubivcek2020relaxation}, modeling microstructure in materials \cite{chipot1988equilibrium, james2005theory, muller2006variational, ball1987fine, bhattacharya2003microstructure}, damage mechanics \cite{rieger2009young, rieger2005model}, optimal design \cite{MR901988, kohn1986optimal,kohn1986optimalII,kohn1986optimalIII,kohn1986thin}, and fluid mechanics \cite{diperna1985compensated, holm1985nonlinear, jordan1995statistical,jordan1997ideal, miller1990statistical, robert1991maximum, lanthaler2015computation}. They arise naturally in characterizing weak limits of sequences of functions and provide a good framework for studying non-convex variational problems. In this section, we aim to review the essential  definitions and their meaning, and discuss how Young measures are useful in optimization problems. 

Borrowed from \cite{ball2005version}, let us first provide an intuitive description of Young measures as a device for characterizing the weak limits of continuous functions. To this end, consider a bounded open set $\Omega \in \R^n$, and a sequence of functions $f^{(n)}: \Omega \to \R$. Choose a point in the domain $x \in \Omega$ and an open set in the co-domain $K \in \R$. For a ball of radius $\delta$ centered at $x$ denoted by $B_\delta(x)$, let us ask the following question: what is the probability that $f^{(n)}$ maps points in $B_\delta(x)$ to set $K$? we can define:

\begin{equation}
    \nu_{x,\delta}^{(n)}(K) = \frac{|\{z \in B_\delta(x): f^{(n)}\in K\}|}{|B_\delta(x)|},
\end{equation}
where $|\cdot|$ is the Lebesgue measure. As $K$ is an arbitrary open set, this defines a probability measure $\nu_{x,\delta}^{(n)}$ on the whole $\R$. Let us first take the limit of $n \to \infty$ and shrink the balls by taking the limit $\delta \to 0$. Under assumptions that will be made precise, we thus have a family of probability measures:
\begin{equation}\label{eq:nu_x}
    \nu_x = \lim_{\delta \to 0} \lim_{n \to \infty} \nu_{x,\delta}^{(n)},
\end{equation}
where the convergence is understood as the weak-$\star$ convergence of probability measures. Thus the Young measure characterizes the local behavior of the sequence of functions. 

To make the above heuristics more precise, we next recall the Young-measure compactness theorem. We denote $C_0(\R)$ as the Banach space of continuous functions vanishing at infinity, and define $\PM(\R)$ as the space of probability measures on $\R$. For a measurable map $x \mapsto \nu_x \in \PM(\R)$, we write $\langle \nu_x, \phi\rangle := \int_\R \phi(z)d\nu_x(z)$.

\begin{theorem}\normalfont{
    Let $(u_j)$ be a bounded sequence in the space of measurable functions $L^1(\Omega ; \R).$ Then there exists a subsequence $(u_{j_k})$ and a measurable family $\nu = (\nu_x)_{x \in \Omega}$ of probability measures on $\R$ such that for every $\phi \in C_0(\R)$ and every $\psi \in L^1(\Omega)$,
    \begin{equation}
        \int_\Omega \psi(x)\phi(u_{j_k}(x))dx \to \int_\Omega \psi(x)\langle\nu_x,\phi\rangle dx.
    \end{equation}
    If additionally $u_{j_k} \rightharpoonup u$ in $L^1(\Omega)$, then the barycenter $\bar{\nu}_x := \langle \nu_x , \text{id}\rangle$ satisfies $\bar{\nu}_x = u(x)$ for almost every $x \in \Omega$.
    }
\end{theorem}\label{thm:compactness}
The proof relies on the Banach--Alaoglu theorem applied to the dual space of \(L^{1}(\Omega; C_0(\mathbb{R}))\) and the separability of \(L^{1}(\Omega; C_0(\mathbb{R}))\).
 We refer the interested reader to \cite{ball2005version, kinderlehrer1991characterizations}. 

Let us now discuss non-convex variational problems and the relevant Young measures. Consider an energy functional $E:X \to \R$, defined as:
\begin{equation}\label{eq:energy}
    E(u) := \int_\Omega W(x, u(x), \nabla u(x))dx,
\end{equation}
where $\Omega$ is a bounded open set in $\R^n$ and $W:\Omega \times \R \times \R^n \to \R$ is the energy density. Suppose $X := \{ u \in W^{1,p}(\Omega) : u=0 \ \text{on} \ \partial \Omega \}$, where $1 < p < \infty$. Consider the variational problem: 
\begin{equation}\label{eq:var_pr}
    \inf_{u \in X} E(u).
\end{equation}
Suppose that the energy density is continuous and can be bounded as:
\begin{equation}\label{eq:Cara}
    a|r|^p-A \le W(x,u,r) \le A(1+|r|^p),
\end{equation}
for some $a, A \in \R$. For many physical applications, this is a reasonable hypothesis which renders $E(u)$ as well-defined, continuous, and coercive. On the other hand, we do not assume $W(x,u, \cdot)$ to be convex. Hence $E$ is not weakly lower semi-continuous. This means that while we know the infimum in \eqref{eq:var_pr} exists, as it is bounded from below due to \eqref{eq:Cara}, the minimizing sequence might become increasingly oscillatory, and do not converge in $X$. Hence, there is no minimizer in the classical sense. 

It is known that an equivalent relaxation can be stated where the minimization is conducted over space of measures \cite{pedregal1997parametrized}. In particular, using gradient Young measure we can write the equivalent problem: 
\begin{equation}\label{eq:var_pr_m}
    \tilde{E}(\nu):= \int_{\Omega}{}\int_{\R^n} W(x,u,\lambda) d\nu_x (\lambda) dx,
\end{equation}
where $\nu \in \tilde{X}$ defined as:
\begin{multline}\label{eq:mspace}
     \tilde{X} = \{ u \in X, \nu: \Omega \to \mathcal{P}(\R^n): \nabla u(x) = \int_{\R^n} \lambda d\nu_x(\lambda), \forall x \in \Omega \\ \ \text{and} \ \int_{\R^n} |\lambda|^p d\nu_x(\lambda) < \infty \}.
\end{multline}
The usefulness of Young measure in the context of non-convex variational problem can be seen from the following known theorem:
\begin{theorem} \normalfont{Suppose $W$ in the expression of energy \eqref{eq:energy} satisfies the Carath\'eodory assumptions \eqref{eq:Cara}. The minimum of $E$, as defined in \eqref{eq:var_pr}, can then be obtained by minimizing $\tilde{E}$ defined in \eqref{eq:var_pr_m}, that is 
    \begin{equation}
   \inf_{u \in X} E(u)= \min_{(u,\nu) \in \tilde{X}} \tilde{E}(\nu).
   \label{eq:relaxation}
\end{equation}}
\end{theorem}

We note that while the original variational form \eqref{eq:var_pr} is not convex, this new formulation \eqref{eq:var_pr_m} is convex with respect to $\nu$. Thus the above theorem is understood as a convex relaxation of the original variational problem. This bonus comes at the cost of changing the search space from $X$ to space of measures $\tilde{X}$ which are harder to approximate. 

We refer the reader interested in the proof to \cite{carstensen2000numerical,ball2005version} or Theorem 4.4 in \cite{pedregal1997parametrized}. To keep the presentation self-contained and clarify key nuances of the subject, we provide a brief sketch of the proof:

\begin{enumerate}
    \item \textit{Minimizing sequence and oscillations:} Since \(W\) is not convex in the gradient variable, minimizing sequences \((u_j) \subset X\) for \(E(u)\) may fail to converge strongly in \(W^{1,p}(\Omega)\) due to oscillations or microstructures developing in \(\nabla u_j\).

    \item \textit{Young measure generation:} By the compactness theorem of Young measures above Thm. \ref{thm:compactness}, up to a subsequence, the gradients \(\nabla u_j\) generate a Young measure \(\nu = (\nu_x)_{x \in \Omega}\), i.e., for every continuous bounded function \(\varphi : \mathbb{R}^n \to \mathbb{R}\),
    \begin{equation}\label{eq:youngmeasureconv}
        \varphi(\nabla u_j(x)) \rightharpoonup \int_{\mathbb{R}^n} \varphi(\lambda) \, d\nu_x(\lambda) \quad \text{weakly in } L^1(\Omega).
    \end{equation}

    \item \textit{Relaxed energy functional:} Using \eqref{eq:youngmeasureconv}, the energy functional along the minimizing sequence satisfies
    \begin{equation}
        \liminf_{j \to \infty} E(u_j) \geq \int_\Omega \int_{\mathbb{R}^n} W(x, u(x), \lambda) \, d\nu_x(\lambda) \, dx =: \tilde{E}(\nu).
    \end{equation}

    \item \textit{Lower semicontinuous envelope and minimization:} The relaxed problem over Young measures \(\nu \in \tilde{X}\) is weakly lower semicontinuous and admits a minimizer, which achieves the same infimum as the original problem:
    \begin{equation}
        \inf_{u \in X} E(u) = \min_{\nu \in \tilde{X}} \tilde{E}(\nu).
    \end{equation}
\end{enumerate}

\section{Neural network approximation of Young measures}

This section proposes a framework for approximating Young measures with a deep neural network (DNN). In a nutshell, the key idea is to represent the unknown Young measure as a push-forward map of a Gaussian measure, and then use DNNs to approximate this push-forward map.

We construct a family of parameterized measure by push-forward as follows. 
Denote a transport map $f: (x,\xi)\in \Omega \times \R^n \mapsto f_x(\xi)\in\R^n$ that characterizes the Young measure as the push forward of a $n$-dimensional standard Gaussian denoted by $\nu_x = (f_x)_{\#} \gamma$. This means that for any Borel set $A \in \R^n$:
\begin{equation} \label{eq:GaussianPushforward}
    \nu_x (A) = \gamma(f_x^{-1}(A)),
\end{equation}
where $f_x^{-1}$ denotes the pre-image of $f_x$ and $\gamma$ is the standard Gaussian distribution on $\R^n$. Denote the density for Gaussian as $\rho(\xi)$, so we have
\begin{equation}\label{eq:pdf}
    \nu_x(A)= \int_{f_x^{-1}(A)} d\gamma(\xi) = \int_{f_x^{-1}(A)} \rho(\xi) d\xi. 
\end{equation}

Using the above construction, we can rewrite the energy functional \eqref{eq:var_pr_m} using the push-forward map as 
\begin{equation}\label{eq:var_pr_m}
    \hat{E}(f):= \int_{\Omega}{}\int_{\R^n} W(x,u,f_x(\xi)) d\gamma(\xi) dx.
\end{equation}

Note that since $f_x$ is the push-forward map for the gradient Young measure \eqref{eq:mspace}, we have to ensure that the expectation of the push-forwarded measure describes a gradient field. To clarify this point, let us define V as follows:
\begin{equation}\label{eq:curlfree}
    V(x) = \int_{\R^n} \lambda\, \nu_x(d\lambda).
\end{equation}
For $V$ to be a gradient, assuming that $\Omega$ is a simply connected domain, $V$ has to be curl-free, i.e.  $\nabla \times V=0$. Leveraging this curl-free condition $u$ can be evaluated at a point $x \in \Omega$ via line integration $u(x) = u(x_0)+\int_0^1 \nabla u(x_0 + t(x-x_0)) \cdot (x-x_0)dt$, where $x_0 \in \partial \Omega$. Hence we can write:
\begin{equation}\label{eq:u}
    u(x) = u(x_0)+\int_0^1 \int_\Omega \int_{\R^n}  f_{x_0 + t(x-x_0)}(\xi)\cdot (x-x_0) d\gamma(\xi)dxdt.
\end{equation}

\medskip 

While we have reformulated the variational problem with respect to parameterized measures to that for a function, the map $f$ is still high dimensional. To proceed, it is natural to consider neural network ansatz for $f$. We first note that $f$ does not have to be continuous. The following one-dimensional example illustrates that:
\begin{example}
    Suppose $\nu = \sum_{i=1}^N a_i \delta_{\xi_i}$, where $\sum_{i=1}^N a_i = 1$. Using \eqref{eq:GaussianPushforward}, one obtains that the forward map $f(\xi)$ must satisfy $f(\xi \in \Omega_i) = \xi_i$, where $\{\Omega_i\}$ with $i \in \{1, 2, ..., N\}$ is a partition of the real line ($\bigcup \Omega_i = \R$) such that $\gamma(\Omega_i)=a_i$.
\end{example}

For the consideration of approximation theory and optimization, it is preferred to parameterize continuous functions, and thus instead of $f$, we use neural network to parameterize a Lipschitz function $F: \Omega \times \R^n \to \R^n$, and take $f = \nabla_{\xi} F$, which is well defined due to the Lipschitzness of $F$.

Consider a deep neural network (DNN) using ResNet architecture \cite{he2016deep} comprised by stacking many blocks where each block consists of two linear transformations, two activation functions, and a residual connection. For an input $x \in \R^{n}$, the $i$-th layer can be expressed as a map $\rho_i$ with $m$ neurons in the residual block:
\begin{equation}\label{eq:resnet}
    \rho_i (x) = \sigma(W_i^2 \cdot \sigma(W_i^1\cdot x+ b^1_i)+b_i^2)+x,
\end{equation}
where $W_i^1\in \R^{m\times n}, W_i^2 \in \R^{n\times m}$ are the weight matrices, and $b_i^1 \in \R^m, b_i^2 \in \R^n$ are the so-called bias vectors. Here $\sigma$ denotes the activation function, which is a nonlinear function that acts on the input vector component-wise, i.e. $\sigma([x_j])=[\sigma(x_j)]$. In essence, if one defines an affine map $a_i^k(x) = W_i^k \cdot x + b_i^k$, then $\rho_i (x) = \sigma \circ a_i^2 \circ \sigma \circ a_i^1 + I$, where $I$ is the identity map. A full ResNet DNN of depth N defines a set of functions:
\begin{equation}\label{eq:DNNRep}
        \mathcal{NN}_\sigma (\theta) = \{ f : f = \rho_N \circ \rho_{N-1} \circ \dots \circ \rho_0 \},
\end{equation}
with $\theta$ representing all the parameters in the network. See Figure~\ref{fig:network_schematic} for an illustration of the architecture. Note that once an activation function is chosen, there are $2N (mn+m+n)$ parameters to be chosen (i.e. weights and biases). The function in such ansatz is Lipschitz as long as $\sigma$ is. 

Let us denote the neural network representation of $F$ as $\hF(x,\xi,\theta)$. By substituting this DNN representation in the minimization problem \eqref{eq:var_pr_m}, we obtain a variational form in terms of the parameters of the DNN:

\begin{equation}\label{eq:var_pr_m_DNN}
  L(\theta):= \int_{\Omega}{}\int_{\R^n} W(x,u,\nabla_{\xi} \hF(x,\xi,\theta)) \, d\gamma(\xi) dx,
\end{equation}
with the optimization problem as $\min_\theta L(\theta)$.

Let us denote the error of approximating Young measure $\nu_x$ using the above procedure as $e(x)$:
\begin{equation}\label{eq:DNN1}
 e(x):=\inf_{\hF(x,\theta) \in \mathcal{NN}_\sigma} \D((\nabla \hF(x,\theta))_{\#} \gamma,\nu_x),
\end{equation}
where $\D(\pi_1,\pi_2)$ is a suitable discrepancy between the probability measures $\pi_1$ and $\pi_2$. The total error $\tilde{e}$ would then be:
\begin{equation}
    \tilde{e} = \int_\Omega e(x).
\end{equation}

Since we are obtaining approximations to  Young measures through variational search using the energy functional, we next show that, for $\D$ chosen as Wasserstein-2 metric, the distance between two Young measures is controlled by their energy.

Let us revisit \eqref{eq:var_pr_m} and denote the minimizer as $\nu^{\star}$: $\tilde{E}[\nu^{\star}] = \inf_{\nu} \tilde{E}[\nu]$. 
Note that 
\begin{equation}
    \begin{aligned}
        \tilde{E}[\nu] - \tilde{E}[\nu^{\star}] & = \int_{\Omega} \int W(x, u, \lambda) - W(x, u, \eta) \, \nu_x(d\lambda)\nu_x^{\star}(d\eta) d x  \\
        &\ge  \int_{\Omega} \inf_{\mu_x \in \Gamma(\nu_x, \nu_x^*)} \int W(x, u, \lambda) - W(x, u, \eta) \, \mu_x(d\lambda, d\eta) d x 
    \end{aligned}
\end{equation}
where for each $x$, $\mu_x$ is optimized over all couplings between  $\nu_x$ and $\nu_x^{\star}$, denoted as $\Gamma(\nu_x, \nu_x^*)$.

If we assume coercivity for stored energy in the sense of:
\begin{equation}
    c \dist(\lambda, \arg\min W(x, u, \cdot))^2 
    \leq W(x, u, \lambda) - \inf W(x, u, \cdot) 
\end{equation}
Then we obtain the desired estimate
\begin{equation}
    \begin{aligned}
    \int_{\Omega} \mathcal{W}_2^2(\nu_x, \nu_x^{\star}) \, dx &=
    \int_{\Omega} \inf_{\mu_x} \int \lVert \lambda - \eta \rVert^2  \, \mu_x(d\lambda, d\eta) d x  \\
    & \leq \frac{1}{c} \int_{\Omega} \inf_{\mu_x} \int W(x, u,\lambda) - W(x, u, \eta)  \, \mu_x(d\lambda, d\eta) d x \\
    & \leq \frac{1}{c} \bigl( \tilde{E}(\nu) - \tilde{E}(\nu^{\star}) \bigr), 
    \end{aligned}
\end{equation}
that is the integrated Wasserstein-$2$ metric between $\nu$ and $\nu^{\star}$ can be controlled by the energy difference. 

On the other hand, thanks to universal approximation theory of neural networks for measures (see e.g., \cite{lu2020universal}), we can guarantee that for a given Young measure, it can be accurately approximated as the width and depth of the DNN ansatz of $F$ become large. Thus the numerical optimization is potentially able to find a good approximate minimizer to the variational problem. We provide more details on the numerical implementation in the next section.

\section{Numerical implementation and  examples}
In this section, we demonstrate several numerical experiments to test the aforementioned framework for approximating Young measures in non-convex variational problems. We first delineate the computational details that are consistently used across all our numerical tests, including parameters of the architecture described in \eqref{eq:resnet} and \eqref{eq:DNNRep}. We then present the numerical results for all our tests in one and two dimensions.\footnote{Here and in the sequel, the dimension means physical dimension; as our variational problems are formulated in terms of Young measures,  they are infinite dimensional.} 

\subsection{Details of neural network implementations}
\label{Sec:PINN}
\begin{figure}[H]
    \centering
    \includegraphics[width=\textwidth]{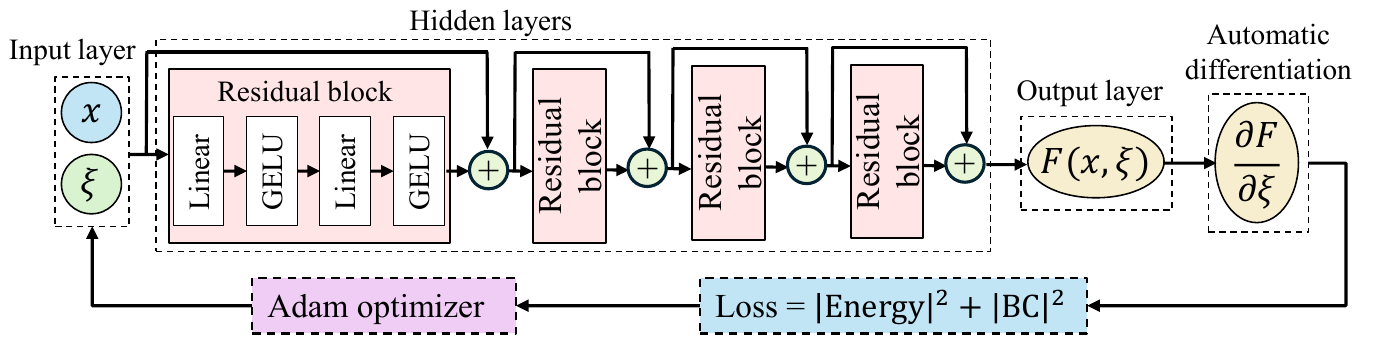} 
    \caption{Schematic of the neural network architecture and computation pipeline.}
    \label{fig:network_schematic}
\end{figure}

The residual architecture (ResNet) architecture \eqref{eq:resnet} \cite{he2016deep} is used with GELU \cite{hendrycks2016gaussian} as the activation function. For a given variational problem, the corresponding ResNet entails a loss function that includes the variational form \eqref{eq:var_pr_m_DNN} and the boundary conditions. 
Hence each loss function includes a domain integral term approximating the energy and penalty terms enforcing boundary conditions, where densities of 1D Gaussian distributions $e^{-\xi^2/2}$, or 2D $e^{-(\xi^2+\tau^2)/2}$, for latent variables $\xi, \tau$ are used in evaluating the integrals. As outlined in the previous section, the neural network learns the map $F$ where $\nabla F$ is the push-forward map from a Gaussian to the Young measure associated to the problem under consideration. Gradients of $F$, such as $\partial F/\partial \xi$ and $\partial F/\partial \tau$ are then computed using automatic differentiation, and these derivatives are directly embedded into the loss function. Figure~\ref{fig:network_schematic} schematically presents the neural network framework, including the input/output layers, residual blocks, and loss computation pipeline.

ResNet is comprised of fully connected residual blocks, which employs skip connections to improve gradient-based training and model expressivity. The network architecture is consistent across all test cases, with minor variations depending on the input dimensionality.

The neural network consists of:
\begin{itemize}
    \item An input layer receiving the coordinates from the \emph{physical domain} and the \emph{Gaussian random variables}. 
    In the one-dimensional case, the inputs are $(x, \xi)$, where $x \in [0,1]$ denotes the physical coordinate and $\xi \sim \mathcal{N}(0,1)$ represents a realization of the Gaussian random variable. 
    For two-dimensional cases the inputs are $(x,y,\xi,\tau)$, where $(x,y)$ denote the physical coordinates and $(\xi,\tau)$ denote Gaussian random variables. 
    \item Four residual blocks, each containing two fully connected layers followed by GELU activation functions and identity skip connections.        
    \item An output layer producing a scalar value $F(x,\xi)$ or $F(x,y,\xi, \tau)$, depending on the dimensionality.
\end{itemize}
Each fully connected layer comprises 25 neurons, and all weights are initialized using Xavier initialization. Model parameters are optimized using the Adam optimizer with an initial learning rate of $10^{-3}$. A ReduceLROnPlateau scheduler \cite{pytorch-reducelronplateau, su2024thermodynamics, liu2025towards} is employed to adaptively lower the learning rate when the loss stagnates. We approximate the loss function by performing importance samplings. Each model is trained over 1000 to 2000 epochs, depending on the experiment. Batch sizes are dynamically adjusted during training in higher-dimensional cases, beginning at 5 and increasing periodically.

\subsection{Case 1: The 1D Bolza problem}

We first investigate the famous Bolza-type example as a canonical non-convex variational problem with known Young measure solutions as:

\begin{equation} \label{eq:1DExample}
\begin{aligned}
    \inf_{u} \quad & E[u] := \int_{0}^1 \bigl((u')^2 -1\bigr)^2 + u^2 dx \\
    \textrm{s.t.} \quad & u(0)=0, \\
    &u(1)=0,
\end{aligned}
\end{equation}

where $u'$ denotes the derivative of $u$. Adopting the framework described in the foregoing section, we can re-write $E$ as:

\begin{equation}
    \hat{E}= \int_0^1 \int_\R \bigl((\frac{dF}{dx}(x,\xi)^2 -1\bigr)^2 + u^2\bigr) d\gamma(\xi)dx.
\end{equation}
We employ the architecture outlined in the previous subsection to develop a neural network approximation of $F$, where we formulate the loss function for the network comprised of the energy term, i.e. the integral in \eqref{eq:1DExample}, and the boundary conditions. We approximate the loss function in by sampling from a uniform grid $x \in [0,1]$ and $\xi \in [-2,2]$ on a $201 \times 201$ grid: 

\begin{equation}\label{eq:1D1}
\begin{aligned}
\text{Loss} = & \underbrace{\frac{1}{N} \sum_{i=1}^N \left(\frac{1}{M} \sum_{k=1}^{M} \left(\left(\frac{dF}{d\xi}(x_i,\xi_k)\right)^2 - 1\right)^2 e^{-\xi_k^2/2}\right)}_{\bigl((u')^2 -1\bigr)^2 \ term \ in \ energy} \\
& + \underbrace{\frac{1}{N} \sum_{i=1}^N \left(\frac{1}{N} \sum_{j=1}^i \frac{1}{M} \sum_{k=1}^M \frac{dF}{d\xi}(x_j,\xi_k) e^{-\xi_k^2/2}\right)^2}_{u^2 \ term \ in \ energy} \\
& + \underbrace{\left(\frac{\lambda}{N} \sum_{i=1}^N \frac{1}{M} \sum_{k=1}^M \frac{dF}{d\xi}(x_i,\xi_k) e^{-\xi_k^2/2}\right)^2}_{right \ boundary \ condition \ u(1)=0}.
\end{aligned}
\end{equation}

\begin{figure}[H]
  \centering
  \begin{subfigure}{0.49\textwidth}
    \includegraphics[width=\linewidth]{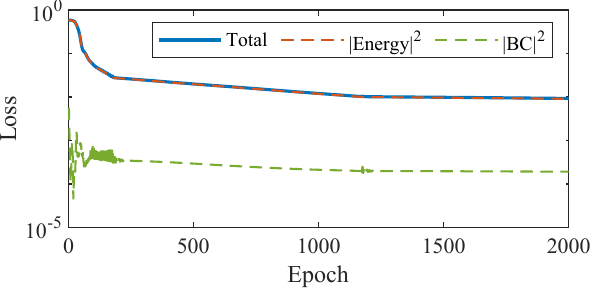}
    \caption{}
    \label{fig:1D1_loss}
  \end{subfigure}\hfill
  \begin{subfigure}{0.49\textwidth}
    \includegraphics[width=\linewidth]{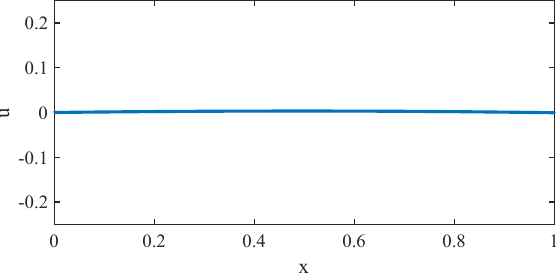}
    \caption{}
    \label{fig:1D1_U}
  \end{subfigure}\hfill
  \begin{subfigure}{0.4\textwidth}
    \includegraphics[width=\linewidth]{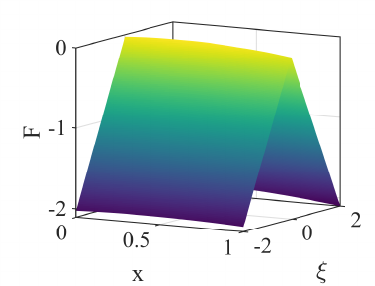}
    \caption{}
    \label{fig:1D1_F3D}
  \end{subfigure}
  \begin{subfigure}{0.4\textwidth}
    \includegraphics[width=\linewidth]{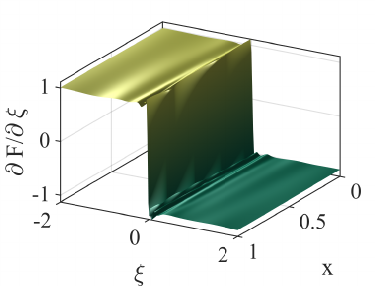}
    \caption{}
    \label{fig:1D1_dFdz3D}
  \end{subfigure}\hfill
  \begin{subfigure}{0.32\textwidth}
    \includegraphics[width=\linewidth]{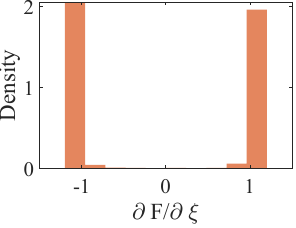}
    \caption{}
    \label{fig:1D1_dFdzhist}
  \end{subfigure}

  \caption{(a) Convergence of total neural network loss during training. (b) Scalar field U obtained via cumulative integration of weighted gradients. (c) 3D surface of \(F\) over \((x, \xi)\).
(d) The push-forward map \(\partial F/\partial \xi\) over \((x, \xi)\).
(e) Distribution of the approximated gradient Young measure.}
  \label{fig:1D1_training}
\end{figure}

The training results are represented in Fig. \ref{fig:1D1_training}, inclduing the We note that using the resultant $F$ we can compute the solution $u$:

\begin{equation}\label{eq:1D1U}
\begin{aligned}
U(x_n) = & \frac{1}{N} \sum_{i=1}^n \left(\frac{1}{M} \sum_{k=1}^{M} \frac{dF}{d\xi}(x_i,\xi_k) e^{-\xi_k^2/2}\right).
\end{aligned}
\end{equation}

In addition to loss over epochs for training, the predicted $F$ and its derivative is illustrated in Fig.~\ref{fig:1D1_F3D} and Fig.~\ref{fig:1D1_dFdz3D}. The histogram in Fig.~\ref{fig:1D1_dFdzhist} represents the distribution of the gradient values $\frac{dF}{d\xi}$. It can be seen that the neural network correctly approximates a homogeneous gradient Young measure for the Bolza problem, with atomic density concentrated at $+1$ and $-1$. 

It can be verified that the computed push-froward measure indeed matches theoretical predictions. We obtain the push forward measure (as shown in Fig. \ref{fig:1D1_dFdz3D}) for every x as:
\[
f_x(\xi)=
\begin{cases}
  1,  & \xi>0,\\
 -1, & \xi<0.
\end{cases}
\]
Using \eqref{eq:pdf}, we then have
\begin{align*}
        &\nu_x(\{1\})=\int_{f_x^{-1}({1})} d\gamma = \frac{1}{\sqrt{2 \pi}}\int_{0}^\infty e^{-\xi^2/2}d\xi = \frac{1}{2}, \\
        &\nu_x(\{-1\})=\int_{f_x^{-1}({-1})} d\gamma = \frac{1}{\sqrt{2 \pi}}\int_{-\infty}^0 e^{-\xi^2/2}d\xi = \frac{1}{2}, 
\end{align*}
which yields homogeneous gradient Young measures for the 1D Bolza problem as
\begin{equation}
    \nu_x = \frac{1}{2}\delta_{\{-1\}} + \frac{1}{2}\delta_{\{1\}}.
\end{equation}

\subsection{Case 2: A quasi one-dimensional problem}

We consider here the following variational problem in a unit square $D=[0,1]^2$:
\begin{equation} \label{eq:Semi2D}
\begin{aligned}
   & \inf_u \int_D \left((u_x^2-1)^2+u_y^2\right)   dA &\\
    \textrm{s.t.} \quad & u(x=0, y)  =0, \quad u(x=1, y)=0,\\
     & u(x, y=0) =0, \quad u(x, y=1)=0.
\end{aligned}
\end{equation}
Note that $u_x$ represents partial derivative with respect to $x$. The competition between non-vanishing partial along x and zero Dirichlet boundary condition leads to emergence of Young measures. Similar to the previous case, we start by re-writing the minimization problem in terms of Young measures that are represented as a push-forward of a Gaussian:
\begin{equation}
    \int_D \int_{\R^2} \bigl( ( \frac{\partial F}{\partial \xi})^2 -1 \bigr)^2 + (\frac{\partial F}{\partial \tau})^2 d\gamma(\xi,\tau) dA.
\end{equation}

Employing the same neural network architecture and training setup detailed in Section~\ref{Sec:PINN}, we set the loss function as:
\begin{equation} \label{eq:LossSemi2D}
\begin{aligned}
\text{Loss} = & \frac{\lambda_1}{N M} \sum_{i=1}^N \sum_{j=1}^M \left\{\frac{1}{R T} \sum_{p=1}^{R} \sum_{q=1}^{T}  \left[ \underbrace{\left(  \left( \frac{\partial F}{\partial \xi}(x_i,y_j, \xi_p, \tau_q) \right)^2 - 1 \right)^2}_{(u_x^2-1)^2 \ term \ in \ energy} \right. \right. \\
& \qquad + \left. \left. \underbrace{\left( \frac{\partial F}{\partial \tau}(x_i,y_j, \xi_p, \tau_q) \right)^2}_{u_y^2 \ term \ in \ energy} \right] e^{-(\xi_p^2+\tau_q^2)/2} \right\}\\
& + \underbrace{{\lambda_2} \sum_{j=1}^M \left( \frac{1}{N} \sum_{i=1}^N \frac{1}{R} \sum_{p=1}^R \frac{1}{T} \sum_{q=1}^T \frac{\partial F}{\partial \xi}(x_i,y_j, \xi_p, \tau_q) e^{-(\xi_p^2+\tau_q^2)/2} \right)^2}_{u(x,1)=0\ boundary \ condition} \\
& + \underbrace{{\lambda_2} \sum_{i=1}^N \left( \frac{1}{M} \sum_{j=1}^M \frac{1}{R} \sum_{p=1}^R \frac{1}{T} \sum_{q=1}^T \frac{\partial F}{\partial \tau}(x_i,y_j, \xi_p, \tau_q) e^{-(\xi_p^2+\tau_q^2)/2} \right)^2}_{u(1,y)=0\ boundary \ condition}\\
& + 
    \frac{\lambda_3}{N} \sum_{i=1}^N \frac{1}{M} \sum_{j=1}^M \left(  \underbrace{\frac{1}{R} \sum_{p=1}^R \frac{1}{T} \sum_{q=1}^T \frac{\partial^2F}{\partial x\, \partial \tau}(x_i,y_j, \xi_p, \tau_q) e^{-(\xi_p^2+\tau_q^2)/2}}_{curl \ free \ condition \ \cdots} \right. \\
&\qquad \left. - \underbrace{\frac{1}{R} \sum_{p=1}^R \frac{1}{T} \sum_{q=1}^T \frac{\partial^2F}{\partial y\, \partial \xi}(x_i,y_j, \xi_p, \tau_q) e^{-(\xi_p^2+\tau_q^2)/2}}_{\cdots \ curl \ free \ condition} \right)^2.
\end{aligned}
\end{equation}

Let us remark that the boundary condition terms in the loss function represent the following realizations of zero Dirichlet condition on the right and top boundary:
\begin{equation} \label{eq:Semi2DBC}
\begin{aligned}
    u(1,y) &= \int_{0}^1 u_x(t,y) \, dt, \\
    u(x,1) &= \int_{0}^1 u_y(x,t) \, dt.
\end{aligned}
\end{equation}
Furthermore, the curl-free term in the loss mirrors conditions described in the previous section and Eq. \eqref{eq:curlfree}. In particular, adopting the same notation as \eqref{eq:curlfree} for this case 2 we have:
\begin{equation*}
    V := [V_1 \, V_2]^T := \int_{\R^2}\lambda\, (\nabla F_{\#} \gamma)( d \lambda).
\end{equation*}
The curl-free condition in the loss function is imposing $\frac{\partial V_2}{\partial x} - \frac{\partial V_1}{\partial y}=0$. We also remark that weights $\lambda_1$, $\lambda_2$, and $\lambda_3$ in the loss function is used as tuning hyperparameters in the minimization process. 

We approximate the loss function in all the 2D cases by performing a Monte Carlo sampling. The \texttt{StochasticMeshgridDataset} Python class is used to generate training batches which randomly samples values from $(x, y, \xi, \tau)$, and  constructs input tensors for efficient parallel loss computation.  

\begin{figure}[H]
  \centering
  \begin{subfigure}[t]{0.55\textwidth}
    \includegraphics[width=\textwidth]{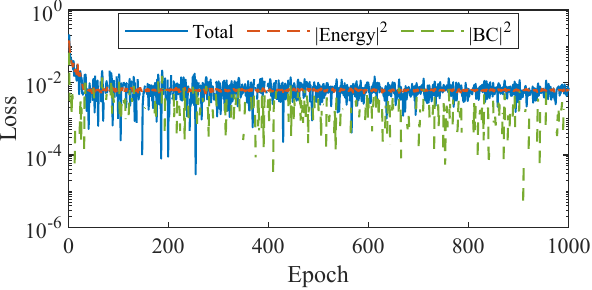}
    \caption{}
    \label{fig:Semi2D_loss}
  \end{subfigure}\hfill
  \begin{subfigure}[t]{0.40\textwidth}
    \includegraphics[width=\textwidth]{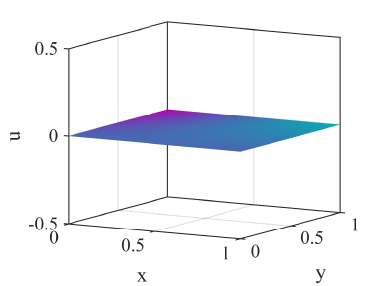}
    \caption{}
    \label{fig:Semi2D_U}
  \end{subfigure}

  \vspace{2ex}

  \begin{subfigure}[b]{0.32\textwidth}
    \includegraphics[width=\textwidth]{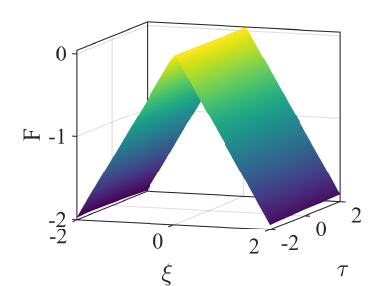}
    \caption{}
    \label{fig:FSemi2D_a}
  \end{subfigure}\hfill
  \begin{subfigure}[b]{0.32\textwidth}
    \includegraphics[width=\textwidth]{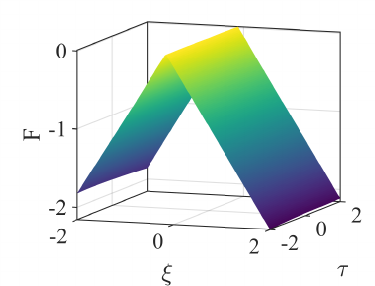}
    \caption{}
    \label{fig:FSemi2D_b}
  \end{subfigure}\hfill
  \begin{subfigure}[b]{0.32\textwidth}
    \includegraphics[width=\textwidth]{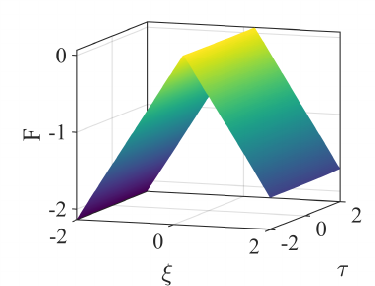}
    \caption{}
    \label{fig:FSemi2D_c}
  \end{subfigure}
  \caption{Top row: (a) Convergence of the neural network’s total loss in log scale.      (b) Predicted scalar field \(u(x,y)\).  
           Bottom row shows the neural network predictions of the field \(F(\xi,\tau)\) at three representative points in the unit square \((x,y)=(0.5,0.5),(0.25,0.75),(0.75,0.25)\), respectively for (c)–(e).}
  \label{fig:Semi2D_merged1}
\end{figure}

\begin{figure}[H]
  \centering
  \begin{subfigure}[b]{0.32\textwidth}
    \includegraphics[width=\textwidth]{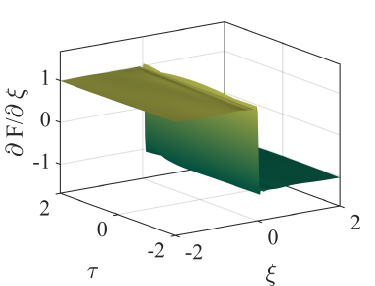}
    \caption{}
    \label{fig:FSemi2D_d}
  \end{subfigure}\hfill
  \begin{subfigure}[b]{0.32\textwidth}
    \includegraphics[width=\textwidth]{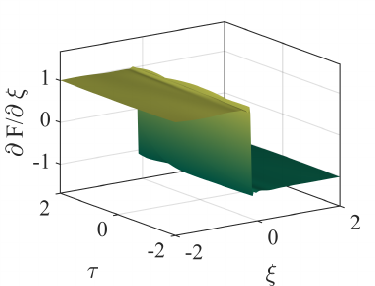}
    \caption{}
    \label{fig:FSemi2D_e}
  \end{subfigure}\hfill
  \begin{subfigure}[b]{0.32\textwidth}
    \includegraphics[width=\textwidth]{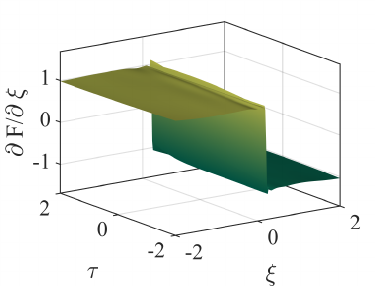}
    \caption{}
    \label{fig:FSemi2D_f}
  \end{subfigure}

  \vspace{1ex}

  \begin{subfigure}[b]{0.32\textwidth}
    \includegraphics[width=\textwidth]{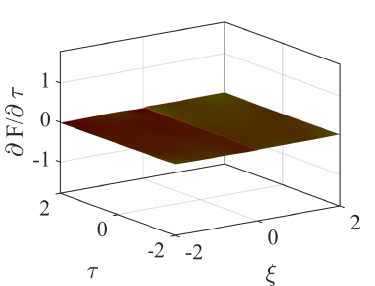}
    \caption{}
    \label{fig:FSemi2D_g}
  \end{subfigure}\hfill
  \begin{subfigure}[b]{0.32\textwidth}
    \includegraphics[width=\textwidth]{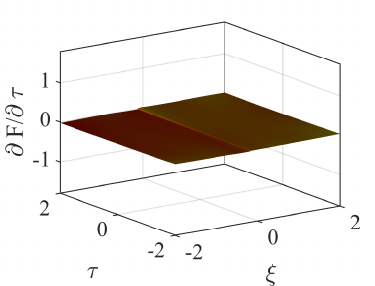}
    \caption{}
    \label{fig:FSemi2D_h}
  \end{subfigure}\hfill
  \begin{subfigure}[b]{0.32\textwidth}
    \includegraphics[width=\textwidth]{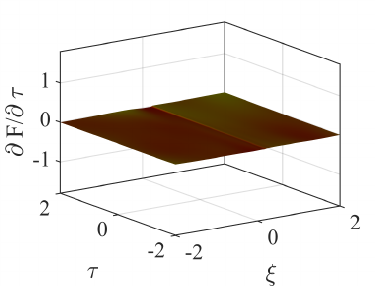}
    \caption{}
    \label{fig:FSemi2D_i}
  \end{subfigure}

  \vspace{2ex}

  \begin{subfigure}[b]{0.32\textwidth}
    \includegraphics[width=\textwidth]{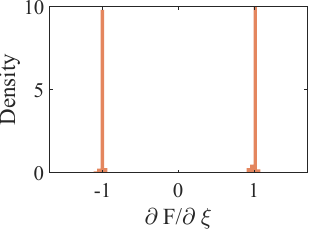}
    \caption{}
    \label{fig:HSemi2D_a}
  \end{subfigure}\hfill
  \begin{subfigure}[b]{0.32\textwidth}
    \includegraphics[width=\textwidth]{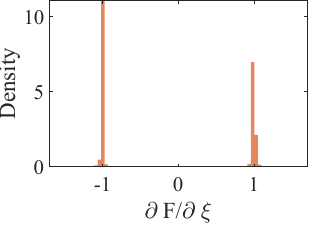}
    \caption{}
    \label{fig:HSemi2D_b}
  \end{subfigure}\hfill
  \begin{subfigure}[b]{0.32\textwidth}
    \includegraphics[width=\textwidth]{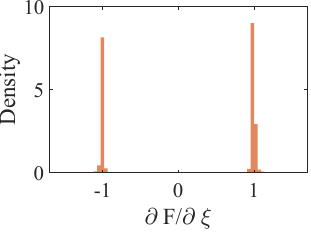}
    \caption{}
    \label{fig:HSemi2D_c}
  \end{subfigure}

  \vspace{1ex}

  \begin{subfigure}[b]{0.32\textwidth}
    \includegraphics[width=\textwidth]{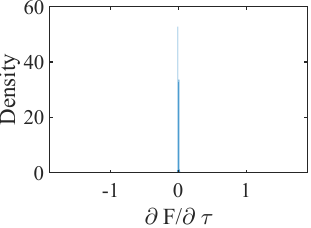}
    \caption{}
    \label{fig:HSemi2D_d}
  \end{subfigure}\hfill
  \begin{subfigure}[b]{0.32\textwidth}
    \includegraphics[width=\textwidth]{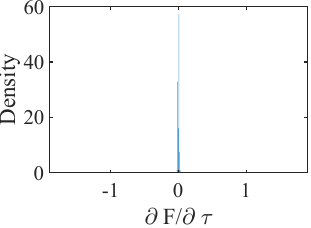}
    \caption{}
    \label{fig:HSemi2D_e}
  \end{subfigure}\hfill
  \begin{subfigure}[b]{0.32\textwidth}
    \includegraphics[width=\textwidth]{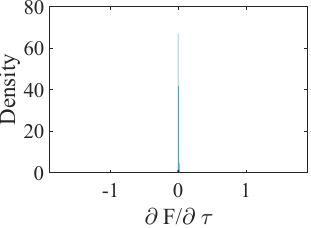}
    \caption{}
    \label{fig:HSemi2D_f}
  \end{subfigure}

  \caption{Top two rows: Components of the two dimensional $\nabla F$ as the push-forward map, depicted respectively at three representative points \((x,y)=(0.5,0.5),(0.25,0.75),(0.75,0.25)\).   
           Bottom two rows: gradient Young measure densities obtained as push-forward of a Gaussian obtained as histograms using 10,000 Gaussian samples.} 
  \label{fig:Semi2D_merged2}
\end{figure}

The results of neural network approximation, including the training loss, $F$, push-forward $\nabla F$, and predicted solution $u$ are shown in Fig.~\ref{fig:Semi2D_merged1} and Fig.~\ref{fig:Semi2D_merged2}.

The results demonstrate that the effective solution is virtually zero everywhere. The gradient Young measures for $u_x$ has atomic mass at $+1$ and $-1$ and the the gradient Young measures for $u_y$ is essentially $\delta_{\{0\}}$. We note that the approximated push-forward map, and the resultant Young measure, is not exactly homogeneous, which we believe can be attributed to the numerical difficulties of computing the four-dimensional $F$. We suspect that further efforts in tuning architecture and hyperparameters of the DNN can lead to smaller loss values and improvements in the push-forward maps. 

The numerical results of case 2 as a quasi-1D problem also matches theoretical expectations, with the approximated Young measure for a fixed $y$ being identical to the 1D Bolza problem in case 1.

\subsection{Case 3: a 2D symmetric four-well problem}

This case considers a variational problem in the unit square $D=[0,1]^2$ defined as:

\begin{equation} \label{eq:2D}
\begin{aligned}
   & \inf_u \int_D \left((u_x^2-1)^2+(u_y^2-1)^2 \right)  dA\\
    \textrm{s.t.} \quad & u(x=0, y)=0, \quad u(x=1, y)=0,\\
    & u(x, y=0)=0, \quad u(x, y=1)=0.
\end{aligned}
\end{equation}

This is a four-well energy landscape, and the minimization of energy term favors u with slopes $\pm 1$, while the boundary conditions enforce a contradicting condition. Combination of these two effects leads again to emergence of parametrized measures. Given the similarities of this case with case 2, we directly write the loss function that is used in the training of the deep neural network for this problem as: 


\begin{equation} \label{eq:loss2D}
\begin{aligned}
\text{Loss} =\ & \frac{\lambda_1}{N M} \sum_{i=1}^N \sum_{j=1}^M \left\{ \frac{1}{R T} \sum_{p=1}^{R} \sum_{q=1}^{T} \left[ \underbrace{\left( \left( \frac{\partial F}{\partial \xi}(x_i,y_j, \xi_p, \tau_q) \right)^2 - 1 \right)^2}_{(u_x^2-1)^2 \ term \ in \ energy}  \right. \right. \\
& \qquad \left. + \underbrace{\left( \left( \frac{\partial F}{\partial \tau}(x_i,y_j, \xi_p, \tau_q) \right)^2 - 1 \right)^2}_{(u_y^2-1)^2 \ term \ in \ energy} \right] e^{-(\xi_p^2+\tau_q^2)/2} \left. \vphantom{\frac{1}{R T}} \right\} \\
& + \underbrace{\lambda_2 \sum_{j=1}^M \left( \frac{1}{N} \sum_{i=1}^N \frac{1}{R} \sum_{p=1}^R \frac{1}{T} \sum_{q=1}^T \frac{\partial F}{\partial \xi}(x_i,y_j, \xi_p, \tau_q) e^{-(\xi_p^2+\tau_q^2)/2} \right)^2}_{u(x,1)=0\ boundary \ condition} \\
& + \underbrace{{\lambda_2}\sum_{i=1}^N \left( \frac{1}{M} \sum_{j=1}^M \frac{1}{R} \sum_{p=1}^R \frac{1}{T} \sum_{q=1}^T \frac{\partial F}{\partial \tau}(x_i,y_j, \xi_p, \tau_q) e^{-(\xi_p^2+\tau_q^2)/2} \right)^2}_{u(1,y)=0\ boundary \ condition}\\
& + 
    \frac{\lambda_3}{N} \sum_{i=1}^N \frac{1}{M} \sum_{j=1}^M \left(  \underbrace{\frac{1}{R} \sum_{p=1}^R \frac{1}{T} \sum_{q=1}^T \frac{\partial^2F}{\partial x\, \partial \tau}(x_i,y_j, \xi_p, \tau_q) e^{-(\xi_p^2+\tau_q^2)/2}}_{curl \ free \ condition \ \cdots} \right. \\
&\qquad \left. - \underbrace{\frac{1}{R} \sum_{p=1}^R \frac{1}{T} \sum_{q=1}^T \frac{\partial^2F}{\partial y\, \partial \xi}(x_i,y_j, \xi_p, \tau_q) e^{-(\xi_p^2+\tau_q^2)/2}}_{\cdots \ curl \ free \ condition} \right)^2.
\end{aligned}
\end{equation}

\begin{figure}[H]
  \centering
  \begin{subfigure}[t]{0.55\textwidth}
    \includegraphics[width=\textwidth]{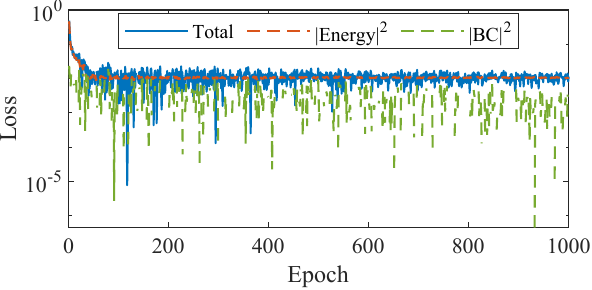}
    \caption{}
    \label{fig:2D3_loss}
  \end{subfigure}\hfill
  \begin{subfigure}[t]{0.40\textwidth}
    \includegraphics[width=\textwidth]{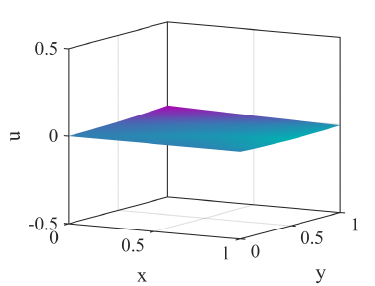}
    \caption{}
    \label{fig:2D3_U}
  \end{subfigure}

  \vspace{2ex}

  \begin{subfigure}[b]{0.32\textwidth}
    \includegraphics[width=\textwidth]{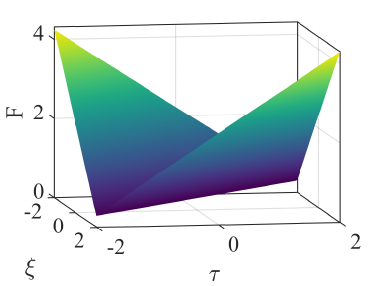}
    \caption{}
    \label{fig:2D3_Fa}
  \end{subfigure}\hfill
  \begin{subfigure}[b]{0.32\textwidth}
    \includegraphics[width=\textwidth]{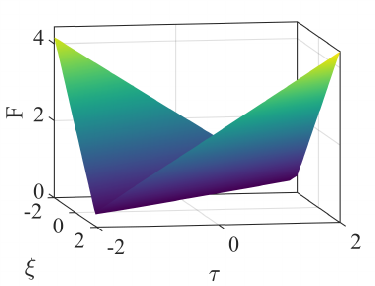}
    \caption{}
    \label{fig:2D3_Fb}
  \end{subfigure}\hfill
  \begin{subfigure}[b]{0.32\textwidth}
    \includegraphics[width=\textwidth]{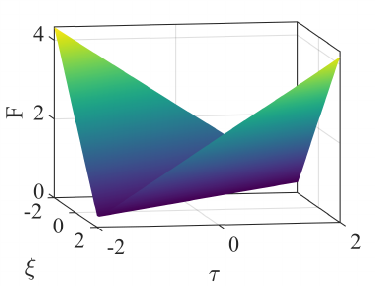}
    \caption{}
    \label{fig:2D3_Fc}
  \end{subfigure}

  \caption{Top row: (a) Convergence of the neural network’s total loss in log scale.  
    (b) Predicted scalar field \(u(x,y)\).  
           Bottom row shows the neural network predictions of the field \(F(\xi,\tau)\) at three representative points in the unit square \((x,y)=(0.5,0.5),(0.25,0.75),(0.75,0.25)\), respectively for (c)–(e).}
  \label{fig:2D3_merged1}
\end{figure}

\begin{figure}[H]
  \centering
  \begin{subfigure}[b]{0.32\textwidth}
    \includegraphics[width=\textwidth]{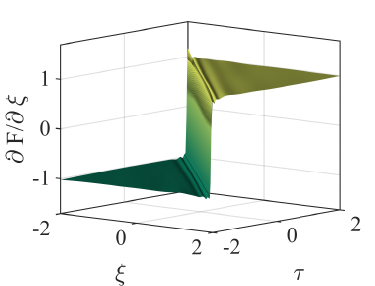}
    \caption{}
    \label{fig:2D3_d}
  \end{subfigure}\hfill
  \begin{subfigure}[b]{0.32\textwidth}
    \includegraphics[width=\textwidth]{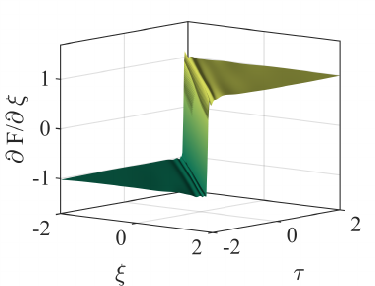}
    \caption{}
    \label{fig:2D3_e}
  \end{subfigure}\hfill
  \begin{subfigure}[b]{0.32\textwidth}
    \includegraphics[width=\textwidth]{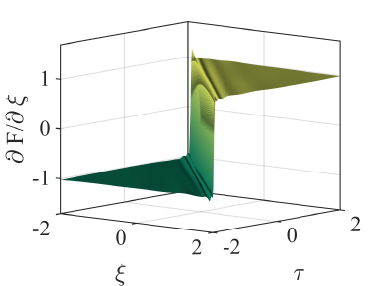}
    \caption{}
    \label{fig:2D3_f}
  \end{subfigure}

  \vspace{1ex}

  \begin{subfigure}[b]{0.32\textwidth}
    \includegraphics[width=\textwidth]{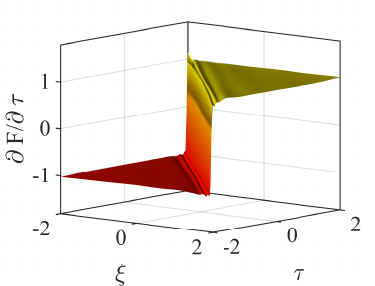}
    \caption{}
    \label{fig:2D3_g}
  \end{subfigure}\hfill
  \begin{subfigure}[b]{0.32\textwidth}
    \includegraphics[width=\textwidth]{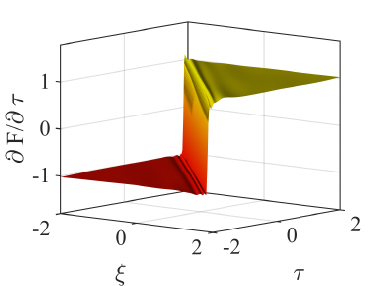}
    \caption{}
    \label{fig:2D3_h}
  \end{subfigure}\hfill
  \begin{subfigure}[b]{0.32\textwidth}
    \includegraphics[width=\textwidth]{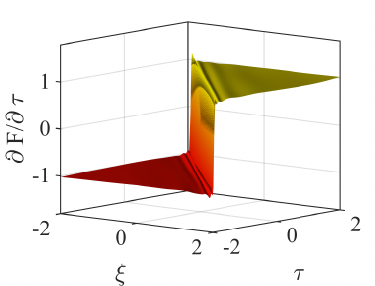}
    \caption{}
    \label{fig:2D3_i}
  \end{subfigure}

  \vspace{2ex}

  \begin{subfigure}[b]{0.32\textwidth}
    \includegraphics[width=\textwidth]{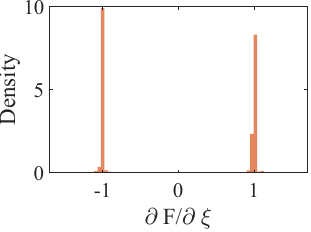}
    \caption{}
    \label{fig:2D3_hist_a}
  \end{subfigure}\hfill
  \begin{subfigure}[b]{0.32\textwidth}
    \includegraphics[width=\textwidth]{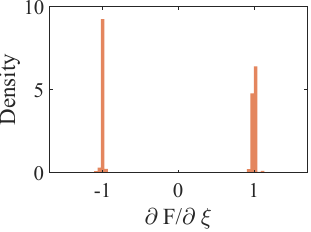}
    \caption{}
    \label{fig:2D3_hist_b}
  \end{subfigure}\hfill
  \begin{subfigure}[b]{0.32\textwidth}
    \includegraphics[width=\textwidth]{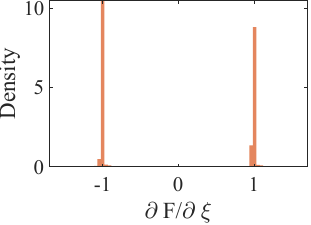}
    \caption{}
    \label{fig:2D3_hist_c}
  \end{subfigure}

  \vspace{1ex}

  \begin{subfigure}[b]{0.32\textwidth}
    \includegraphics[width=\textwidth]{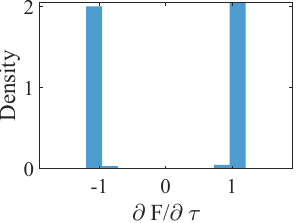}
    \caption{}
    \label{fig:2D3_hist_d}
  \end{subfigure}\hfill
  \begin{subfigure}[b]{0.32\textwidth}
    \includegraphics[width=\textwidth]{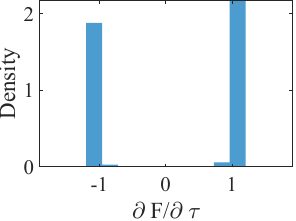}
    \caption{}
    \label{fig:2D3_hist_e}
  \end{subfigure}\hfill
  \begin{subfigure}[b]{0.32\textwidth}
    \includegraphics[width=\textwidth]{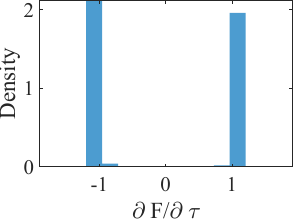}
    \caption{}
    \label{fig:2D3_hist_f}
  \end{subfigure}

  \caption{Top two rows: Components of the two dimensional $\nabla F$ as the push-forward map, depicted respectively at three representative points \((x,y)=(0.5,0.5),(0.25,0.75),(0.75,0.25)\).   
           Bottom two rows: gradient Young measure densities obtained as push-forward of a Gaussian obtained as histograms using 10,000 Gaussian samples.}
  \label{fig:2D3_merged2}
\end{figure}

In Fig.~\ref{fig:2D3_merged1} and Fig.~\ref{fig:2D3_merged2}, we show the training loss and neural network solution to Eq.~\ref{eq:loss2D} along with the approximated $F$, partial derivatives of $F$, predicted solution $u$, and densities for the corresponding Young measure.

It can be seen that the gradient Young measure describing the minimizing sequence of the variational problem \eqref{eq:2D} is approximated as Dirac distributions with support at partial derivatives equal to $+1$ and $-1$. We again remark that further efforts can be carried out to reduce the training loss via different stochastic gradient descent algorithms or further tuning of hyperparameters such as learning rate. We also believe that such improvements will not have a major effect in the results reported in Fig.~\ref{fig:2D3_merged1} and Fig.~\ref{fig:2D3_merged2}.

\subsection{Case 4: a 2D non-symmetric two-well problem}
This final numerical experiment focuses on a variational problem that case where neither the effective solution nor the underlying Young measure is known \textit{a priori}. In this case, we restrict ourselves again to a unit square $D=[0,1]^2$ and define the following variational problem as:  
\begin{equation} \label{eq:2DSinintegral}
\begin{aligned}
    &\inf_u \int_D \left((u_x^2-1)^2+u_y^2\right)   dA \\
    \textrm{s.t.} \quad & u(x=0, y) =0, \quad u(x=1, y)=\alpha y,\\ &u(x, y=0)=0, \quad u(x, y=1)=\alpha x,
\end{aligned}
\end{equation}
where we choose $\alpha = 10^{-2}$. The boundary conditions oppose the bulk energy terms that favor a solution that is flat along the y-direction. This incompatibility between the \textit{bulk-boundary constraints} again leads to emergence of Young measures. 

While the loss function for the corresponding neural network resembles case 2, except for the top and right boundary condition, we write the full loss function for the sake of completeness:
\begin{equation} \label{eq:loss2DSin}
\begin{aligned}
\text{Loss} = & \frac{\lambda_1}{N M} \sum_{i=1}^N \sum_{j=1}^M \left\{\frac{1}{R T} \sum_{p=1}^{R} \sum_{q=1}^{T}  \left[ \underbrace{\left(  \left( \frac{\partial F}{\partial \xi}(x_i,y_j, \xi_p, \tau_q) \right)^2 - 1 \right)^2}_{(u_x^2-1)^2 \ term \ in \ energy} \right. \right. \\
& \qquad + \left. \left. \underbrace{\left( \frac{\partial F}{\partial \tau}(x_i,y_j, \xi_p, \tau_q) \right)^2}_{u_y^2 \ term \ in \ energy} \right] e^{-(\xi_p^2+\tau_q^2)/2} \right\}\\
& + \lambda_2 \underbrace{\sum_{j=1}^M \left( \frac{1}{N} \sum_{i=1}^N \frac{1}{R} \sum_{p=1}^R \frac{1}{T} \sum_{q=1}^T \frac{\partial F}{\partial \xi}(x_i,y_j, \xi_p, \tau_q) e^{-(\xi_p^2+\tau_q^2)/2 } -\alpha y_j\right)^2}_{u(1,y)=\alpha y\ boundary \ condition} \\
& + \lambda_2 \underbrace{\sum_{i=1}^N \left( \frac{1}{M} \sum_{j=1}^M \frac{1}{R} \sum_{p=1}^R \frac{1}{T} \sum_{q=1}^T \frac{\partial F}{\partial \tau}(x_i,y_j, \xi_p, \tau_q) e^{-(\xi_p^2+\tau_q^2)/2} -\alpha x_i\right)^2}_{u(x,1)=\alpha x\ boundary \ condition}\\
& + 
    \frac{\lambda_3}{N} \sum_{i=1}^N \frac{1}{M} \sum_{j=1}^M \left(  \underbrace{\frac{1}{R} \sum_{p=1}^R \frac{1}{T} \sum_{q=1}^T \frac{\partial^2F}{\partial x\, \partial \tau}(x_i,y_j, \xi_p, \tau_q) e^{-(\xi_p^2+\tau_q^2)/2}}_{curl \ free \ condition \ \cdots} \right. \\
& \qquad \left. - \underbrace{\frac{1}{R} \sum_{p=1}^R \frac{1}{T} \sum_{q=1}^T \frac{\partial^2F}{\partial y\, \partial \xi}(x_i,y_j, \xi_p, \tau_q) e^{-(\xi_p^2+\tau_q^2)/2}}_{\cdots \ curl \ free \ condition} \right)^2.
\end{aligned}
\end{equation}

\begin{figure}[H]
  \centering
  \begin{subfigure}[t]{0.55\textwidth}
    \includegraphics[width=\textwidth]{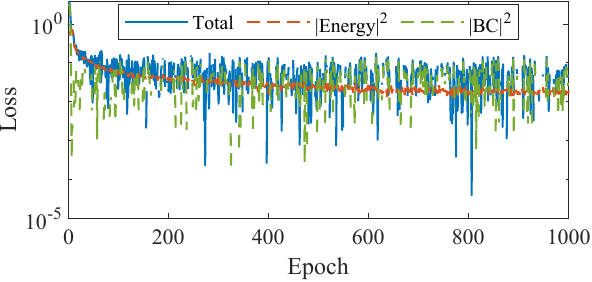}
    \caption{}
    \label{fig:2DSin_loss}
  \end{subfigure}\hfill
  \begin{subfigure}[t]{0.40\textwidth}
    \includegraphics[width=\textwidth]{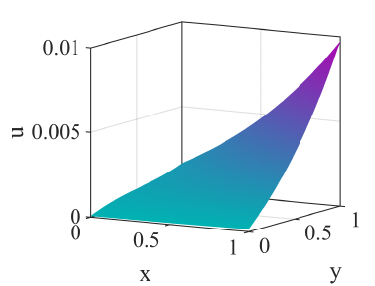}
    \caption{}
    \label{fig:2DSin_U}
  \end{subfigure}

  \vspace{2ex}

  \begin{subfigure}[b]{0.32\textwidth}
    \includegraphics[width=\textwidth]{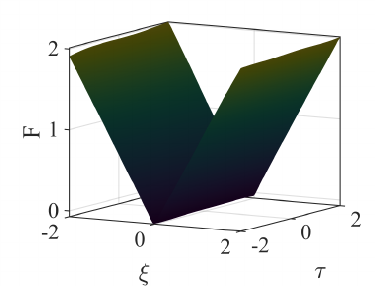}
    \caption{}
    \label{fig:2DSin_Fa}
  \end{subfigure}\hfill
  \begin{subfigure}[b]{0.32\textwidth}
    \includegraphics[width=\textwidth]{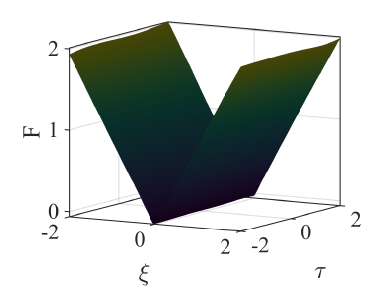}
    \caption{}
    \label{fig:2DSin_Fb}
  \end{subfigure}\hfill
  \begin{subfigure}[b]{0.32\textwidth}
    \includegraphics[width=\textwidth]{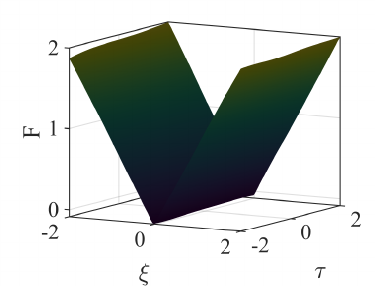}
    \caption{}
    \label{fig:2DSin_Fc}
  \end{subfigure}

  \caption{Top row: (a) Convergence of the neural network’s total loss in log scale.  
    (b) Predicted scalar field \(u(x,y)\).  
           Bottom row shows the neural network predictions of the field \(F(\xi,\tau)\) at three representative points in the unit square \((x,y)=(0.5,0.5),(0.25,0.75),(0.75,0.25)\), respectively for (c)–(e).}
  \label{fig:2DSin_merged1}
\end{figure}

\begin{figure}[H]
  \centering
  \begin{subfigure}[b]{0.32\textwidth}
    \includegraphics[width=\textwidth]{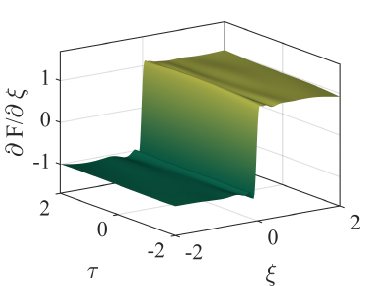}
    \caption{}
    \label{fig:2DSin_d}
  \end{subfigure}\hfill
  \begin{subfigure}[b]{0.32\textwidth}
    \includegraphics[width=\textwidth]{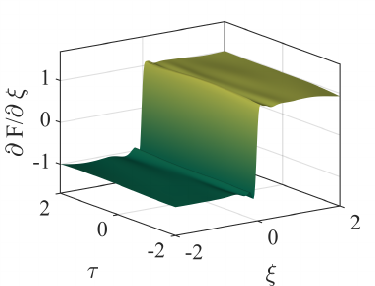}
    \caption{}
    \label{fig:2DSin_e}
  \end{subfigure}\hfill
  \begin{subfigure}[b]{0.32\textwidth}
    \includegraphics[width=\textwidth]{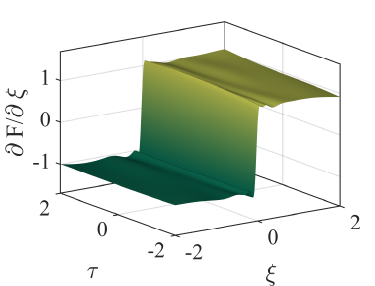}
    \caption{}
    \label{fig:2DSin_f}
  \end{subfigure}

  \vspace{1ex}

  \begin{subfigure}[b]{0.32\textwidth}
    \includegraphics[width=\textwidth]{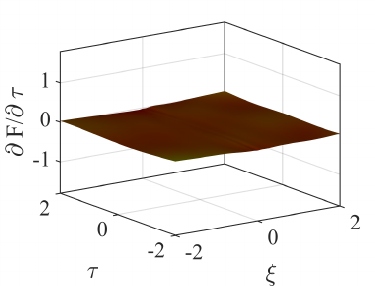}
    \caption{}
    \label{fig:2DSin_g}
  \end{subfigure}\hfill
  \begin{subfigure}[b]{0.32\textwidth}
    \includegraphics[width=\textwidth]{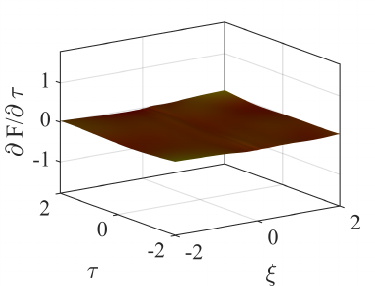}
    \caption{}
    \label{fig:2DSin_h}
  \end{subfigure}\hfill
  \begin{subfigure}[b]{0.32\textwidth}
    \includegraphics[width=\textwidth]{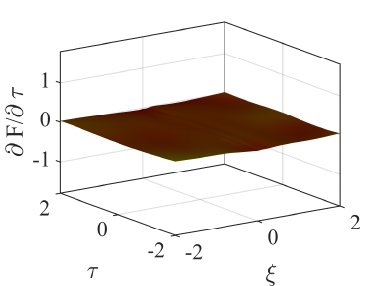}
    \caption{}
    \label{fig:2DSin_i}
  \end{subfigure}

  \vspace{2ex}

  \begin{subfigure}[b]{0.32\textwidth}
    \includegraphics[width=\textwidth]{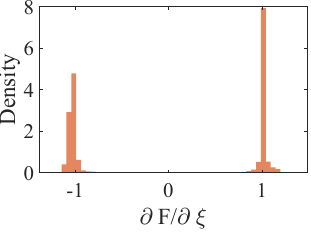}
    \caption{}
    \label{fig:2DSin_hist_a}
  \end{subfigure}\hfill
  \begin{subfigure}[b]{0.32\textwidth}
    \includegraphics[width=\textwidth]{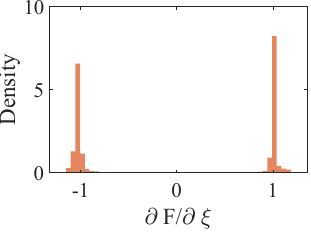}
    \caption{}
    \label{fig:2DSin_hist_b}
  \end{subfigure}\hfill
  \begin{subfigure}[b]{0.32\textwidth}
    \includegraphics[width=\textwidth]{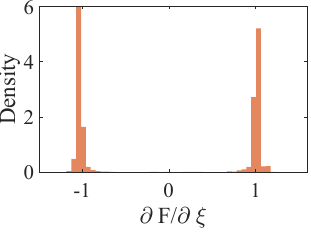}
    \caption{}
    \label{fig:2DSin_hist_c}
  \end{subfigure}

  \vspace{1ex}

  \begin{subfigure}[b]{0.32\textwidth}
    \includegraphics[width=\textwidth]{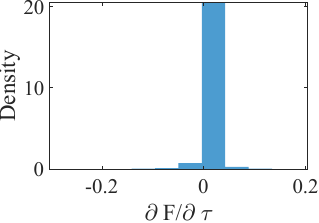}
    \caption{}
    \label{fig:2DSin_hist_d}
  \end{subfigure}\hfill
  \begin{subfigure}[b]{0.32\textwidth}
    \includegraphics[width=\textwidth]{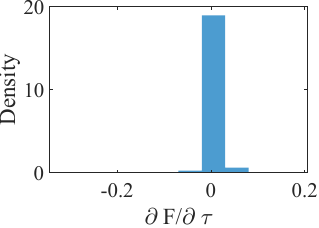}
    \caption{}
    \label{fig:2DSin_hist_e}
  \end{subfigure}\hfill
  \begin{subfigure}[b]{0.32\textwidth}
    \includegraphics[width=\textwidth]{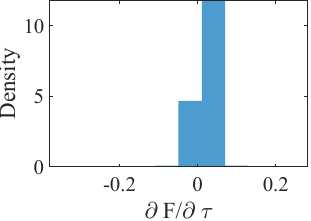}
    \caption{}
    \label{fig:2DSin_hist_f}
  \end{subfigure}

  \caption{Top two rows: Components of the two dimensional $\nabla F$ as the push-forward map, depicted respectively at three representative points \((x,y)=(0.5,0.5),(0.25,0.75),(0.75,0.25)\).   
           Bottom two rows: gradient Young measure densities obtained as push-forward of a Gaussian obtained as histograms using 10,000 Gaussian samples.}
  \label{fig:2DSin_merged2}
\end{figure}

Fig.~\ref{fig:2DSin_merged1} and Fig.~\ref{fig:2DSin_merged2} summarize the results for case 4. We find out that the computed u has values very close to zero. We also observe that the distribution of $u_x$ has concentrations on $\pm1$.

Although no closed-form Young measure solution is available for this case, the learned measure reproduces key qualitative features: in the bulk, the distribution of \(u_x\) is bimodal with mass near \(\{\pm 1\}\), and the field \(u\) remains close to zero across most of \(D\). By contrast, the learned distribution of \(u_y\) concentrates near $0$ across the domain.  

Based on our experiments, further optimizer and hyperparameter sweeps primarily shift training loss values without altering the robust traits (bimodality in \(u_x\), near-vanishing \(u\) in the interior).

\section{Conclusion and Outlook}

A powerful idea in the direct method of the calculus of variations is to enlarge the class of admissible minimizing sequences from functions to measures, thereby
replacing explicit \textit{(quasi)-convexification} of the integrand by the emergence of parameterized (Young) measures generated by minimizing sequences. Motivated by
the ability of deep neural networks (DNNs) to approximate high–dimensional objects, we proposed a neural representation of Young measures in this paper: each
$\nu_x$ is modeled as the pushforward of a simple base law (here, a Gaussian) through a learned transport map $f_\theta(x,\cdot)$, so that
$\nu_x = (f_\theta(x,\cdot))_\#\mathcal N(0,I)$. Combined with the classical observation that finite Gaussian mixtures are weakly dense in $\mathcal P(\mathbb R^d)$, this pushforward parameterization yields a practical
and expressive scheme for approximating Young measures in nonconvex problems.

After reviewing the relevant theory in Section~2, we detailed the proposed construction and training objectives in Section~3, including enforcement of barycentric (gradient–Young) admissibility and physics-informed penalties. We then applied the framework to four nonconvex variational problems. Progressing in difficulty, we began with a 1D Bolza problem and proceeded to two-dimensional
settings; our final example tackles a case with no \emph{a priori} known Young-measure solution. Across these experiments, we demonstrated that the Young measures are are directly approximated and the effective solutions are produced without going through the widely-used relaxation or convexification route.

We envision that this framework can be readily extended to the vector-valued fields $u:\Omega\to\mathbb{R}^m$, where the relevant Young measures live on $\mathbb{R}^{m\times d}$, as the natural next step. Such an extension would enable data-driven modeling of microstructure in multi-well elastic energies—including martensitic phase transformations \cite{bhattacharya2003microstructure}—by
learning mixtures over variant wells, recovering
laminate hierarchies and volume fractions directly from the learned measures. We hope this paves the way for a new program for modeling and discovering microstructure in materials such as alloys \cite{song2013enhanced}. 

\section*{Data availability}
The code and data can be found at this GitHub repository: \url{https://github.com/RayeheKM/DNN-YoungMeasure}.

\section*{Acknowledgment}

The work of R.KM. and H.S. are in part supported by Duke University internal Beyond the Horizon grant. The work of J.L.~is supported in part by National Science Foundation via award DMS-2309378.

\bibliographystyle{elsarticle-num-names} 
\bibliography{refs}

\begin{thebibliography}{57}
\expandafter\ifx\csname natexlab\endcsname\relax\def\natexlab#1{#1}\fi
\providecommand{\url}[1]{\texttt{#1}}
\providecommand{\href}[2]{#2}
\providecommand{\path}[1]{#1}
\providecommand{\DOIprefix}{doi:}
\providecommand{\ArXivprefix}{arXiv:}
\providecommand{\URLprefix}{URL: }
\providecommand{\Pubmedprefix}{pmid:}
\providecommand{\doi}[1]{\href{http://dx.doi.org/#1}{\path{#1}}}
\providecommand{\Pubmed}[1]{\href{pmid:#1}{\path{#1}}}
\providecommand{\bibinfo}[2]{#2}
\ifx\xfnm\relax \def\xfnm[#1]{\unskip,\space#1}\fi
\bibitem[{Goodman et~al.(1986)Goodman, Kohn, and Reyna}]{goodman1986numerical}
\bibinfo{author}{J.~Goodman}, \bibinfo{author}{R.~V. Kohn},
  \bibinfo{author}{L.~Reyna},
\newblock \bibinfo{title}{Numerical study of a relaxed variational problem from
  optimal design},
\newblock \bibinfo{journal}{Computer Methods in Applied Mechanics and
  Engineering} \bibinfo{volume}{57} (\bibinfo{year}{1986})
  \bibinfo{pages}{107--127}.
\bibitem[{De~Simone(1993)}]{de1993energy}
\bibinfo{author}{A.~De~Simone},
\newblock \bibinfo{title}{Energy minimizers for large ferromagnetic bodies},
\newblock \bibinfo{journal}{Archive for rational mechanics and analysis}
  \bibinfo{volume}{125} (\bibinfo{year}{1993}) \bibinfo{pages}{99--143}.
\bibitem[{Luskin(1996)}]{luskin1996computation}
\bibinfo{author}{M.~Luskin},
\newblock \bibinfo{title}{On the computation of crystalline microstructure},
\newblock \bibinfo{journal}{Acta numerica} \bibinfo{volume}{5}
  (\bibinfo{year}{1996}) \bibinfo{pages}{191--257}.
\bibitem[{Ball and James(1987)}]{ball1987fine}
\bibinfo{author}{J.~M. Ball}, \bibinfo{author}{R.~D. James},
\newblock \bibinfo{title}{Fine phase mixtures as minimizers of energy},
\newblock \bibinfo{journal}{Archive for Rational Mechanics and Analysis}
  \bibinfo{volume}{100} (\bibinfo{year}{1987}) \bibinfo{pages}{13--52}.
\bibitem[{Chipot and Kinderlehrer(1988)}]{chipot1988equilibrium}
\bibinfo{author}{M.~Chipot}, \bibinfo{author}{D.~Kinderlehrer},
\newblock \bibinfo{title}{Equilibrium configurations of crystals},
\newblock \bibinfo{journal}{Archive for Rational Mechanics and Analysis}
  \bibinfo{volume}{103} (\bibinfo{year}{1988}) \bibinfo{pages}{237--277}.
\bibitem[{Bhattacharya(2003)}]{bhattacharya2003microstructure}
\bibinfo{author}{K.~Bhattacharya}, \bibinfo{title}{Microstructure of
  martensite: why it forms and how it gives rise to the shape-memory effect},
  volume~\bibinfo{volume}{2}, \bibinfo{publisher}{Oxford University Press},
  \bibinfo{year}{2003}.
\bibitem[{Kohn(1991)}]{kohn1991relaxation}
\bibinfo{author}{R.~V. Kohn},
\newblock \bibinfo{title}{The relaxation of a double-well energy},
\newblock \bibinfo{journal}{Continuum Mechanics and Thermodynamics}
  \bibinfo{volume}{3} (\bibinfo{year}{1991}) \bibinfo{pages}{193--236}.
\bibitem[{Govindjee and Miehe(2001)}]{govindjee2001multi}
\bibinfo{author}{S.~Govindjee}, \bibinfo{author}{C.~Miehe},
\newblock \bibinfo{title}{A multi-variant martensitic phase transformation
  model: formulation and numerical implementation},
\newblock \bibinfo{journal}{Computer Methods in Applied Mechanics and
  Engineering} \bibinfo{volume}{191} (\bibinfo{year}{2001})
  \bibinfo{pages}{215--238}.
\bibitem[{Ortiz and Repetto(1999)}]{ortiz1999nonconvex}
\bibinfo{author}{M.~Ortiz}, \bibinfo{author}{E.~Repetto},
\newblock \bibinfo{title}{Nonconvex energy minimization and dislocation
  structures in ductile single crystals},
\newblock \bibinfo{journal}{Journal of the Mechanics and Physics of Solids}
  \bibinfo{volume}{47} (\bibinfo{year}{1999}) \bibinfo{pages}{397--462}.
\bibitem[{Aubry et~al.(2003)Aubry, Fago, and Ortiz}]{aubry2003constrained}
\bibinfo{author}{S.~Aubry}, \bibinfo{author}{M.~Fago},
  \bibinfo{author}{M.~Ortiz},
\newblock \bibinfo{title}{A constrained sequential-lamination algorithm for the
  simulation of sub-grid microstructure in martensitic materials},
\newblock \bibinfo{journal}{Computer Methods in Applied Mechanics and
  Engineering} \bibinfo{volume}{192} (\bibinfo{year}{2003})
  \bibinfo{pages}{2823--2843}.
\bibitem[{Nicolaides and Walkington(1993)}]{nicolaides1993computation}
\bibinfo{author}{R.~A. Nicolaides}, \bibinfo{author}{N.~J. Walkington},
\newblock \bibinfo{title}{Computation of microstructure utilizing young measure
  representations},
\newblock \bibinfo{journal}{Journal of intelligent material systems and
  structures} \bibinfo{volume}{4} (\bibinfo{year}{1993})
  \bibinfo{pages}{457--462}.
\bibitem[{Carstensen and Roub{\'\i}{\v{c}}ek(2000)}]{carstensen2000numerical}
\bibinfo{author}{C.~Carstensen}, \bibinfo{author}{T.~Roub{\'\i}{\v{c}}ek},
\newblock \bibinfo{title}{Numerical approximation of young measuresin
  non-convex variational problems},
\newblock \bibinfo{journal}{Numerische Mathematik} \bibinfo{volume}{84}
  (\bibinfo{year}{2000}) \bibinfo{pages}{395--415}.
\bibitem[{Raissi et~al.(2019)Raissi, Perdikaris, and
  Karniadakis}]{raissi2019physics}
\bibinfo{author}{M.~Raissi}, \bibinfo{author}{P.~Perdikaris},
  \bibinfo{author}{G.~E. Karniadakis},
\newblock \bibinfo{title}{Physics-informed neural networks: A deep learning
  framework for solving forward and inverse problems involving nonlinear
  partial differential equations},
\newblock \bibinfo{journal}{Journal of Computational physics}
  \bibinfo{volume}{378} (\bibinfo{year}{2019}) \bibinfo{pages}{686--707}.
\bibitem[{Karniadakis et~al.(2021)Karniadakis, Kevrekidis, Lu, Perdikaris,
  Wang, and Yang}]{karniadakis2021physics}
\bibinfo{author}{G.~E. Karniadakis}, \bibinfo{author}{I.~G. Kevrekidis},
  \bibinfo{author}{L.~Lu}, \bibinfo{author}{P.~Perdikaris},
  \bibinfo{author}{S.~Wang}, \bibinfo{author}{L.~Yang},
\newblock \bibinfo{title}{Physics-informed machine learning},
\newblock \bibinfo{journal}{Nature Reviews Physics} \bibinfo{volume}{3}
  (\bibinfo{year}{2021}) \bibinfo{pages}{422--440}.
\bibitem[{Han et~al.(2018)Han, Jentzen, and E}]{han2018solving}
\bibinfo{author}{J.~Han}, \bibinfo{author}{A.~Jentzen}, \bibinfo{author}{W.~E},
\newblock \bibinfo{title}{Solving high-dimensional partial differential
  equations using deep learning},
\newblock \bibinfo{journal}{Proceedings of the National Academy of Sciences}
  \bibinfo{volume}{115} (\bibinfo{year}{2018}) \bibinfo{pages}{8505--8510}.
\bibitem[{Beck et~al.(2019)Beck, E, and Jentzen}]{beck2019machine}
\bibinfo{author}{C.~Beck}, \bibinfo{author}{W.~E},
  \bibinfo{author}{A.~Jentzen},
\newblock \bibinfo{title}{Machine learning approximation algorithms for
  high-dimensional fully nonlinear partial differential equations and
  second-order backward stochastic differential equations},
\newblock \bibinfo{journal}{Journal of Nonlinear Science} \bibinfo{volume}{29}
  (\bibinfo{year}{2019}) \bibinfo{pages}{1563--1619}.
\bibitem[{Li et~al.(2020)Li, Kovachki, Azizzadenesheli, Liu, Bhattacharya,
  Stuart, and Anandkumar}]{li2020fourier}
\bibinfo{author}{Z.~Li}, \bibinfo{author}{N.~Kovachki},
  \bibinfo{author}{K.~Azizzadenesheli}, \bibinfo{author}{B.~Liu},
  \bibinfo{author}{K.~Bhattacharya}, \bibinfo{author}{A.~Stuart},
  \bibinfo{author}{A.~Anandkumar},
\newblock \bibinfo{title}{Fourier neural operator for parametric partial
  differential equations},
\newblock \bibinfo{journal}{arXiv preprint arXiv:2010.08895}
  (\bibinfo{year}{2020}).
\bibitem[{Kovachki et~al.(2023)Kovachki, Li, Liu, Azizzadenesheli,
  Bhattacharya, Stuart, and Anandkumar}]{kovachki2023neural}
\bibinfo{author}{N.~Kovachki}, \bibinfo{author}{Z.~Li},
  \bibinfo{author}{B.~Liu}, \bibinfo{author}{K.~Azizzadenesheli},
  \bibinfo{author}{K.~Bhattacharya}, \bibinfo{author}{A.~Stuart},
  \bibinfo{author}{A.~Anandkumar},
\newblock \bibinfo{title}{Neural operator: Learning maps between function
  spaces with applications to pdes},
\newblock \bibinfo{journal}{Journal of Machine Learning Research}
  \bibinfo{volume}{24} (\bibinfo{year}{2023}) \bibinfo{pages}{1--97}.
\bibitem[{Khoo et~al.(2021)Khoo, Lu, and Ying}]{khoo2021solving}
\bibinfo{author}{Y.~Khoo}, \bibinfo{author}{J.~Lu}, \bibinfo{author}{L.~Ying},
\newblock \bibinfo{title}{Solving parametric pde problems with artificial
  neural networks},
\newblock \bibinfo{journal}{European Journal of Applied Mathematics}
  \bibinfo{volume}{32} (\bibinfo{year}{2021}) \bibinfo{pages}{421--435}.
\bibitem[{Lu and Lu(2020)}]{lu2020universal}
\bibinfo{author}{Y.~Lu}, \bibinfo{author}{J.~Lu},
\newblock \bibinfo{title}{A universal approximation theorem of deep neural
  networks for expressing probability distributions},
\newblock \bibinfo{journal}{Advances in neural information processing systems}
  \bibinfo{volume}{33} (\bibinfo{year}{2020}) \bibinfo{pages}{3094--3105}.
\bibitem[{Young(1937)}]{young1937generalized}
\bibinfo{author}{L.~C. Young},
\newblock \bibinfo{title}{Generalized curves and the existence of an attained
  absolute minimum in the calculus of variations},
\newblock \bibinfo{journal}{Comptes Rendus de la Societe des Sci. et des
  Lettres de Varsovie} \bibinfo{volume}{30} (\bibinfo{year}{1937})
  \bibinfo{pages}{212--234}.
\bibitem[{Young(2024)}]{young2024lectures}
\bibinfo{author}{L.~C. Young}, \bibinfo{title}{Lectures on the calculus of
  variations and optimal control theory}, volume \bibinfo{volume}{304},
  \bibinfo{publisher}{American Mathematical Society}, \bibinfo{year}{2024}.
\bibitem[{Tartar(1983)}]{tartar1983compensated}
\bibinfo{author}{L.~Tartar},
\newblock \bibinfo{title}{The compensated compactness method applied to systems
  of conservation laws},
\newblock in: \bibinfo{booktitle}{Systems of nonlinear partial differential
  equations}, \bibinfo{publisher}{Springer}, \bibinfo{year}{1983}, pp.
  \bibinfo{pages}{263--285}.
\bibitem[{Tartar(1990)}]{tartar1990h}
\bibinfo{author}{L.~Tartar},
\newblock \bibinfo{title}{H-measures, a new approach for studying
  homogenisation, oscillations and concentration effects in partial
  differential equations},
\newblock \bibinfo{journal}{Proceedings of the Royal Society of Edinburgh
  Section A: Mathematics} \bibinfo{volume}{115} (\bibinfo{year}{1990})
  \bibinfo{pages}{193--230}.
\bibitem[{Allaire and Kohn(1994)}]{allaire1994optimal}
\bibinfo{author}{G.~Allaire}, \bibinfo{author}{R.~V. Kohn},
\newblock \bibinfo{title}{Optimal lower bounds on the elastic energy of a
  composite made from two non-well-ordered isotropic materials},
\newblock \bibinfo{journal}{Quarterly of applied mathematics}
  \bibinfo{volume}{52} (\bibinfo{year}{1994}) \bibinfo{pages}{311--333}.
\bibitem[{Ball(1989)}]{ball2005version}
\bibinfo{author}{J.~M. Ball},
\newblock \bibinfo{title}{A version of the fundamental theorem for {Y}oung
  measures},
\newblock in: \bibinfo{editor}{J.~M. Ball} (Ed.), \bibinfo{booktitle}{PDEs and
  Continuum Models of Phase Transitions}, volume \bibinfo{volume}{344} of
  \textit{\bibinfo{series}{Lecture Notes in Physics}},
  \bibinfo{publisher}{Springer}, \bibinfo{year}{1989}, pp.
  \bibinfo{pages}{207--215}. \URLprefix
  \url{https://doi.org/10.1007/BFb0075577}. \DOIprefix\doi{10.1007/BFb0075577}.
\bibitem[{Francfort and Milton(1994)}]{francfort1994sets}
\bibinfo{author}{G.~A. Francfort}, \bibinfo{author}{G.~W. Milton},
\newblock \bibinfo{title}{Sets of conductivity and elasticity tensors stable
  under lamination},
\newblock \bibinfo{journal}{Communications on pure and applied mathematics}
  \bibinfo{volume}{47} (\bibinfo{year}{1994}) \bibinfo{pages}{257--279}.
\bibitem[{Milton(1990)}]{milton1990characterizing}
\bibinfo{author}{G.~W. Milton},
\newblock \bibinfo{title}{On characterizing the set of possible effective
  tensors of composites: the variational method and the translation method},
\newblock \bibinfo{journal}{Communications on Pure and Applied Mathematics}
  \bibinfo{volume}{43} (\bibinfo{year}{1990}) \bibinfo{pages}{63--125}.
\bibitem[{Milton(1994)}]{milton1994link}
\bibinfo{author}{G.~W. Milton},
\newblock \bibinfo{title}{A link between sets of tensors stable under
  lamination and quasiconvexity},
\newblock \bibinfo{journal}{Communications on Pure and Applied Mathematics}
  \bibinfo{volume}{47} (\bibinfo{year}{1994}) \bibinfo{pages}{959--1003}.
\bibitem[{Nesi and Milton(1991)}]{nesi1991polycrystalline}
\bibinfo{author}{V.~Nesi}, \bibinfo{author}{G.~W. Milton},
\newblock \bibinfo{title}{Polycrystalline configurations that maximize
  electrical resistivity},
\newblock \bibinfo{journal}{Journal of the Mechanics and Physics of Solids}
  \bibinfo{volume}{39} (\bibinfo{year}{1991}) \bibinfo{pages}{525--542}.
\bibitem[{Alberti and M{\"u}ller(2001)}]{alberti2001new}
\bibinfo{author}{G.~Alberti}, \bibinfo{author}{S.~M{\"u}ller},
\newblock \bibinfo{title}{A new approach to variational problems with multiple
  scales},
\newblock \bibinfo{journal}{Communications on Pure and Applied Mathematics: A
  Journal Issued by the Courant Institute of Mathematical Sciences}
  \bibinfo{volume}{54} (\bibinfo{year}{2001}) \bibinfo{pages}{761--825}.
\bibitem[{Balder(1984)}]{balder1984general}
\bibinfo{author}{E.~J. Balder},
\newblock \bibinfo{title}{A general approach to lower semicontinuity and lower
  closure in optimal control theory},
\newblock \bibinfo{journal}{SIAM journal on control and optimization}
  \bibinfo{volume}{22} (\bibinfo{year}{1984}) \bibinfo{pages}{570--598}.
\bibitem[{Roub{\'\i}{\v{c}}ek(2020)}]{roubivcek2020relaxation}
\bibinfo{author}{T.~Roub{\'\i}{\v{c}}ek}, \bibinfo{title}{Relaxation in
  optimization theory and variational calculus}, volume~\bibinfo{volume}{4},
  \bibinfo{publisher}{Walter de Gruyter GmbH \& Co KG}, \bibinfo{year}{2020}.
\bibitem[{James and Kinderlehrer(2005)}]{james2005theory}
\bibinfo{author}{R.~James}, \bibinfo{author}{D.~Kinderlehrer},
\newblock \bibinfo{title}{Theory of diffusionless phase transitions},
\newblock in: \bibinfo{booktitle}{PDEs and Continuum Models of Phase
  Transitions: Proceedings of an NSF-CNRS Joint Seminar Held in Nice, France,
  January 18--22, 1988}, \bibinfo{organization}{Springer},
  \bibinfo{year}{2005}, pp. \bibinfo{pages}{51--84}.
\bibitem[{M{\"u}ller(2006)}]{muller2006variational}
\bibinfo{author}{S.~M{\"u}ller},
\newblock \bibinfo{title}{Variational models for microstructure and phase
  transitions},
\newblock in: \bibinfo{booktitle}{Calculus of Variations and Geometric
  Evolution Problems: Lectures given at the 2nd Session of the Centro
  Internazionale Matematico Estivo (CIME) held in Cetraro, Italy, June 15--22,
  1996}, \bibinfo{publisher}{Springer}, \bibinfo{year}{2006}, pp.
  \bibinfo{pages}{85--210}.
\bibitem[{Rieger and Zimmer(2009)}]{rieger2009young}
\bibinfo{author}{M.~O. Rieger}, \bibinfo{author}{J.~Zimmer},
\newblock \bibinfo{title}{Young measure flow as a model for damage},
\newblock \bibinfo{journal}{Zeitschrift f{\"u}r angewandte Mathematik und
  Physik} \bibinfo{volume}{60} (\bibinfo{year}{2009}) \bibinfo{pages}{1--32}.
\bibitem[{Rieger(2005)}]{rieger2005model}
\bibinfo{author}{M.~O. Rieger},
\newblock \bibinfo{title}{A model for hysteresis in mechanics using local
  minimizers of young measures},
\newblock in: \bibinfo{booktitle}{Elliptic and Parabolic Problems: A Special
  Tribute to the Work of Haim Brezis}, \bibinfo{publisher}{Springer},
  \bibinfo{year}{2005}, pp. \bibinfo{pages}{403--414}.
\bibitem[{Bonnetier and Vogelius(1987)}]{MR901988}
\bibinfo{author}{E.~Bonnetier}, \bibinfo{author}{M.~Vogelius},
\newblock \bibinfo{title}{Relaxation of a compliance functional for a plate
  optimization problem},
\newblock in: \bibinfo{booktitle}{Applications of multiple scaling in mechanics
  ({P}aris, 1986)}, volume~\bibinfo{volume}{4} of
  \textit{\bibinfo{series}{Rech. Math. Appl.}}, \bibinfo{publisher}{Masson,
  Paris}, \bibinfo{year}{1987}, pp. \bibinfo{pages}{31--53}.
\bibitem[{Kohn and Strang(1986{\natexlab{a}})}]{kohn1986optimal}
\bibinfo{author}{R.~V. Kohn}, \bibinfo{author}{G.~Strang},
\newblock \bibinfo{title}{Optimal design and relaxation of variational
  problems, i},
\newblock \bibinfo{journal}{Communications on pure and applied mathematics}
  \bibinfo{volume}{39} (\bibinfo{year}{1986}{\natexlab{a}})
  \bibinfo{pages}{113--137}.
\bibitem[{Kohn and Strang(1986{\natexlab{b}})}]{kohn1986optimalII}
\bibinfo{author}{R.~V. Kohn}, \bibinfo{author}{G.~Strang},
\newblock \bibinfo{title}{Optimal design and relaxation of variational
  problems, ii},
\newblock \bibinfo{journal}{Communications on Pure and Applied Mathematics}
  \bibinfo{volume}{39} (\bibinfo{year}{1986}{\natexlab{b}})
  \bibinfo{pages}{139--182}.
\bibitem[{Kohn and Strang(1986{\natexlab{c}})}]{kohn1986optimalIII}
\bibinfo{author}{R.~V. Kohn}, \bibinfo{author}{G.~Strang},
\newblock \bibinfo{title}{Optimal design and relaxation of variational
  problems, iii},
\newblock \bibinfo{journal}{Communications on Pure and Applied Mathematics}
  \bibinfo{volume}{39} (\bibinfo{year}{1986}{\natexlab{c}})
  \bibinfo{pages}{353--377}.
\bibitem[{Kohn and Vogelius(1986)}]{kohn1986thin}
\bibinfo{author}{R.~V. Kohn}, \bibinfo{author}{M.~Vogelius},
\newblock \bibinfo{title}{Thin plates with rapidly varying thickness, and their
  relation to structural optimization},
\newblock in: \bibinfo{booktitle}{Homogenization and effective moduli of
  materials and media}, \bibinfo{publisher}{Springer}, \bibinfo{year}{1986},
  pp. \bibinfo{pages}{126--149}.
\bibitem[{DiPerna(1985)}]{diperna1985compensated}
\bibinfo{author}{R.~J. DiPerna},
\newblock \bibinfo{title}{Compensated compactness and general systems of
  conservation laws},
\newblock \bibinfo{journal}{Transactions of the American mathematical society}
  \bibinfo{volume}{292} (\bibinfo{year}{1985}) \bibinfo{pages}{383--420}.
\bibitem[{Holm et~al.(1985)Holm, Marsden, Ratiu, and
  Weinstein}]{holm1985nonlinear}
\bibinfo{author}{D.~D. Holm}, \bibinfo{author}{J.~E. Marsden},
  \bibinfo{author}{T.~Ratiu}, \bibinfo{author}{A.~Weinstein},
\newblock \bibinfo{title}{Nonlinear stability of fluid and plasma equilibria},
\newblock \bibinfo{journal}{Physics reports} \bibinfo{volume}{123}
  (\bibinfo{year}{1985}) \bibinfo{pages}{1--116}.
\bibitem[{Jordan(1995)}]{jordan1995statistical}
\bibinfo{author}{R.~Jordan},
\newblock \bibinfo{title}{A statistical equilibrium model of coherent
  structures in magnetohydrodynamics},
\newblock \bibinfo{journal}{Nonlinearity} \bibinfo{volume}{8}
  (\bibinfo{year}{1995}) \bibinfo{pages}{585}.
\bibitem[{Jordan and Turkington(1997)}]{jordan1997ideal}
\bibinfo{author}{R.~Jordan}, \bibinfo{author}{B.~Turkington},
\newblock \bibinfo{title}{Ideal magnetofluid turbulence in two dimensions},
\newblock \bibinfo{journal}{Journal of statistical physics}
  \bibinfo{volume}{87} (\bibinfo{year}{1997}) \bibinfo{pages}{661--695}.
\bibitem[{Miller(1990)}]{miller1990statistical}
\bibinfo{author}{J.~Miller},
\newblock \bibinfo{title}{Statistical mechanics of euler equations in two
  dimensions},
\newblock \bibinfo{journal}{Physical review letters} \bibinfo{volume}{65}
  (\bibinfo{year}{1990}) \bibinfo{pages}{2137}.
\bibitem[{Robert(1991)}]{robert1991maximum}
\bibinfo{author}{R.~Robert},
\newblock \bibinfo{title}{A maximum-entropy principle for two-dimensional
  perfect fluid dynamics},
\newblock \bibinfo{journal}{Journal of Statistical Physics}
  \bibinfo{volume}{65} (\bibinfo{year}{1991}) \bibinfo{pages}{531--553}.
\bibitem[{Lanthaler and Mishra(2015)}]{lanthaler2015computation}
\bibinfo{author}{S.~Lanthaler}, \bibinfo{author}{S.~Mishra},
\newblock \bibinfo{title}{Computation of measure-valued solutions for the
  incompressible euler equations},
\newblock \bibinfo{journal}{Mathematical Models and Methods in Applied
  Sciences} \bibinfo{volume}{25} (\bibinfo{year}{2015})
  \bibinfo{pages}{2043--2088}.
\bibitem[{Kinderlehrer and Pedregal(1991)}]{kinderlehrer1991characterizations}
\bibinfo{author}{D.~Kinderlehrer}, \bibinfo{author}{P.~Pedregal},
\newblock \bibinfo{title}{Characterizations of young measures generated by
  gradients},
\newblock \bibinfo{journal}{Archive for rational mechanics and analysis}
  \bibinfo{volume}{115} (\bibinfo{year}{1991}) \bibinfo{pages}{329--365}.
\bibitem[{Pedregal(1997)}]{pedregal1997parametrized}
\bibinfo{author}{P.~Pedregal}, \bibinfo{title}{Parametrized measures and
  variational principles}, \bibinfo{publisher}{Springer Science \& Business
  Media}, \bibinfo{year}{1997}.
\bibitem[{He et~al.(2016)He, Zhang, Ren, and Sun}]{he2016deep}
\bibinfo{author}{K.~He}, \bibinfo{author}{X.~Zhang}, \bibinfo{author}{S.~Ren},
  \bibinfo{author}{J.~Sun},
\newblock \bibinfo{title}{Deep residual learning for image recognition},
\newblock in: \bibinfo{booktitle}{Proceedings of the IEEE conference on
  computer vision and pattern recognition}, \bibinfo{year}{2016}, pp.
  \bibinfo{pages}{770--778}.
\bibitem[{Hendrycks and Gimpel(2016)}]{hendrycks2016gaussian}
\bibinfo{author}{D.~Hendrycks}, \bibinfo{author}{K.~Gimpel},
\newblock \bibinfo{title}{Gaussian error linear units (gelus)},
\newblock \bibinfo{journal}{arXiv preprint arXiv:1606.08415}
  (\bibinfo{year}{2016}).
\bibitem[{{PyTorch Contributors}(????)}]{pytorch-reducelronplateau}
\bibinfo{author}{{PyTorch Contributors}},
  \bibinfo{title}{torch.optim.lr\_scheduler.ReduceLROnPlateau}, ???? \URLprefix
  \url{https://pytorch.org/docs/stable/generated/torch.optim.lr_scheduler.ReduceLROnPlateau.html},
  \bibinfo{note}{accessed: 2025-08-30}.
\bibitem[{Su et~al.(2024)Su, Yu, Chen, Guo, and Yang}]{su2024thermodynamics}
\bibinfo{author}{M.~Su}, \bibinfo{author}{Y.~Yu}, \bibinfo{author}{T.~Chen},
  \bibinfo{author}{N.~Guo}, \bibinfo{author}{Z.~Yang},
\newblock \bibinfo{title}{A thermodynamics-informed neural network for
  elastoplastic constitutive modeling of granular materials},
\newblock \bibinfo{journal}{Computer Methods in Applied Mechanics and
  Engineering} \bibinfo{volume}{430} (\bibinfo{year}{2024})
  \bibinfo{pages}{117246}.
\bibitem[{Liu et~al.(2025)Liu, Koric, Abueidda, Meidani, and
  Geubelle}]{liu2025towards}
\bibinfo{author}{Q.~Liu}, \bibinfo{author}{S.~Koric},
  \bibinfo{author}{D.~Abueidda}, \bibinfo{author}{H.~Meidani},
  \bibinfo{author}{P.~Geubelle},
\newblock \bibinfo{title}{Towards signed distance function based metamaterial
  design: Neural operator transformer for forward prediction and diffusion
  model for inverse design},
\newblock \bibinfo{journal}{arXiv preprint arXiv:2504.01195}
  (\bibinfo{year}{2025}).
\bibitem[{Song et~al.(2013)Song, Chen, Dabade, Shield, and
  James}]{song2013enhanced}
\bibinfo{author}{Y.~Song}, \bibinfo{author}{X.~Chen},
  \bibinfo{author}{V.~Dabade}, \bibinfo{author}{T.~W. Shield},
  \bibinfo{author}{R.~D. James},
\newblock \bibinfo{title}{Enhanced reversibility and unusual microstructure of
  a phase-transforming material},
\newblock \bibinfo{journal}{Nature} \bibinfo{volume}{502}
  (\bibinfo{year}{2013}) \bibinfo{pages}{85--88}.

\end{thebibliography}

\end{document}